\documentclass[11pt]{article}
\usepackage{amsmath,amssymb,latexsym,float,epsfig,subfigure}
\usepackage{latexsym}
\usepackage{natbib}
\usepackage{amsthm}
\usepackage{cases}
\usepackage{epsfig}
\usepackage{subfigure}
\usepackage{scalefnt}
\usepackage{booktabs}
\usepackage[english]{babel}
\usepackage{graphicx}
\usepackage{titlesec}
\usepackage{epsfig}
\usepackage{enumerate}
\usepackage{caption}
\usepackage{subfigure}
\usepackage{booktabs}
\usepackage[english]{babel}
\usepackage{graphicx}
\usepackage{multirow}
\usepackage{rotating}
\usepackage{fancyhdr}
\usepackage{setspace}
\usepackage[section]{placeins}
\usepackage{xcolor}

\usepackage{breqn}

\usepackage{array, makecell}
\usepackage{pseudocode}
\usepackage{algorithmicx}
\usepackage{algorithm}
\usepackage{url} 
\usepackage{appendix}
\usepackage{ragged2e}
\usepackage{physics}

\usepackage{chngpage}
\newtheorem{theorem}{Theorem}
\newtheorem{proposition}{Proposition}
\newtheorem{definition}{Definition}
\newtheorem{assumption}[theorem]{Assumption}
\newtheorem*{lemma}{Lemma}

\usepackage[nocomma]{optidef}
\usepackage{mathtools,zref-savepos}
\newtagform{roman}[]()

\topmargin by -\headsep \textheight 9in \oddsidemargin 0pt
\evensidemargin \oddsidemargin \marginparwidth 0.5in \textwidth
6.5in
\begin{document}
\title{Queueing Systems with Preferred Service Delivery Times and Multiple Customer Classes}

\author{Melis Boran \thanks{melis.boran@metu.edu.tr}
\and Bahar \c{C}avdar \thanks{bcavdar@tamu.edu} 
\and Tu\u{g}\c{c}e I\c{s}{\i}k \thanks{tisik@clemson.edu}}

\date{$^*$\small{Industrial Engineering Department, Middle East Technical University, Ankara, Turkey }\\
 $^*$\small{Engineering Technology and Industrial Distribution, Texas A\&M University, College Station, TX, USA }\\
    $^\ddag$Department of Industrial Engineering, Clemson University, Clemson, SC, USA\\[2ex]%
}

\maketitle
\begin{abstract}
{
Motivated by the operational problems in click and collect systems, such as curbside pickup programs, we study a joint admission control and capacity allocation problem. We consider a system where arriving customers have preferred service delivery times and gauge the service quality based on the service provider's ability to complete the service as close as possible to the preferred time. Customers can be of different priority classes, and their priority may increase as they wait longer in the queue. The service provider can reject customers upon their arrival if the system is overloaded or outsource the service (alternatively work overtime) when the capacity is not enough. The service provider's goal is to find the minimum-cost admission and capacity allocation policy to dynamically decide when to serve and whom to serve. We model this problem as a Markov Decision Process. Our structural results partially characterize a set of suboptimal solutions, and we develop solution methods using these results. We also develop a problem-specific approximation method that is based on state aggregation to overcome the computational challenges. We present extensive computational results and discuss the impact of problem parameters on the optimal policy.
}
\end{abstract}

\textbf{Keywords:} Capacity allocation, Markov Decision Process, Customer preferences, Curbside pickup, State aggregation, Order fulfillment

\section{Introduction}

Advancements in online platforms and changes in customer expectations have been transforming the market operations and strategies in various industries. In recent years, the traditional retail market has been challenged by the availability of online shopping as digitization redefines what it means to ``go shopping". The shift towards online shopping is not new. In 2015, Nielsen (\cite{Nielsen2015}) reports that about one-fourth of global consumers shop for groceries online, and 55\% say they would buy groceries online. However, the COVID-19 pandemic further accelerated the digitization of the consumer. Today, 7 out of 10 consumers shop online, and 86\% say they would stick to this new habit (\cite{Somers2020}). 
Further, millennials and generation Z, the future drivers of the economy, make about 60\% of their purchases online (\cite{Curbside}). Noticing the upcoming opportunities created by online channels, giant retailers have increasingly been adopting e-commerce, i.e., omnichannel models, and they have become a serious threat to traditional brick-and-mortar retailers such as grocers. Therefore, the traditional brick-and-mortar retailers should adapt their systems, enhance the shopping experience they provide, and meet customers' evolving habits to maintain their competitiveness and market position (\cite{Chopra2016}).

One strategy that has been attracting the attention of giant retailers is establishing ``click and collect" services, which integrates online and offline commerce. These services allow customers to order online for a pickup at a store or another location. Even though the pre-pandemic adoption rate of these services (e.g., in-store, drive-through, and curbside pickup) was just over 10\%, in post-pandemic statistics, the adoption rate increases to 25\% (\cite{Sheehan2021}). Among these service alternatives, curbside pickup, which allows customers to pickup their orders without leaving their car, has already shown its efficiency evolving as a convenience factor for many shoppers (\cite{Curbside}) with its adoption rate rising from 16 percent pre-pandemic to 36 percent in June 2021 (\cite{Standish2020}). Powered by a fast expansion from giant retailers, grocery curbside pickup has become a \$35 billion industry by 2020 (\cite{Forbes}). With more than 2,000 curbside pickup locations, US retailer Walmart leads the market, and the curbside pickup accounts for 33\% of its total digital sales. Kroger, Target, and Whole Foods are in the race to offer curbside pickup services to their customers (\cite{Forbes}).  

Even though mostly giant retailers seem to have succeeded in rolling out curbside pickup services; brick-and-mortar retailers of any type can benefit from developing their curbside pickup capabilities. Curbside pickup seamlessly combines the convenience of online shopping with the trusted physical experience of shopping in-store. These systems have the potential to benefit both retailers and customers. Retailers can provide the convenience, ease, and flexibility of online ordering using their existing resources. In a curbside pickup system, each customer request is received with a preferred service delivery time; hence, customers can plan in advance, order online, and pickup at their convenience much faster compared to last-mile delivery. However, these mutual benefits come at the expense of increased operational challenges to be managed by the retailer, which might make maintaining service quality difficult. An important measure of quality in these systems is the ability to accommodate the customers' preferred order pickup time. Therefore, retailers are forced to make the best use of their resources and plan their operations in a dynamic environment to ensure that customers can be served on time at a low cost. 

Our study is mainly motivated by these challenges that arise in curbside pickup systems. We consider a service system where arriving customers are presented with a fixed-length order horizon to place their order. Upon arrival, customers reveal their preferred service completion time within this order horizon. Customers may belong to different priority classes, such as  high-priority customers who are frequent shoppers or has a paid membership and the remaining lower priority customers. We assume that low-priority customers gain priority as they wait in the system since the retailer aims to provide an acceptable service experience for all customers. From the retailer's perspective, any deviation from the customer's preferred service completion time incurs some cost. In case of an early service completion, the items ordered may need to be stored until pickup or may lose their freshness, resulting in early service costs. In case of a capacity shortage, the retailer can outsource some jobs or increase the capacity temporarily at additional cost. The costs incurred may depend on the customer type. The retailer can also reject low-priority orders based on the current workload and expectations of the future demand. Service provider's problem is to determine the simultaneous admission and capacity allocation decisions considering random arrivals, as well as the associated costs of rejection, early service, and overtime.
We model the this system as a discrete-time Markov Decision Process (MDP), in which we address admission control and the dynamic capacity allocation decisions jointly to minimize long-run average costs. 

Although we mainly consider the retail industry, similar considerations are relevant for other service systems, such as healthcare or maintenance services. In these systems, patients or customers from various groups are served using limited resources, e.g., a diagnostic device. Hence, the service provider must first determine whether to accept or reject a request based on the service capacity. For instance, a private laboratory can receive samples from independent practitioners, insurance companies, clinical research sites, and other health clinics for analysis. Depending on the contract between the laboratory and these parties, the tests may have different priority levels; hence some of the requests may be rejected due to limited equipment capacity. In such systems, there can be a gap between the time the service is requested and the time it is needed. 
Our results also provide useful insights for these systems. 

Here we provide a summary of our key findings. Our results show that the service provider's optimal strategy strongly depends on its service capacity. For small systems that allow customers to place orders over a short time period, the set of suboptimal policies can be easily characterized based on system parameters, such as server capacity and operations costs. Elimination of suboptimal policies significantly increases computational efficiency of traditional solution methods; and the optimal policy can be found within seconds. For large instances, when the quality of the solution is more critical than the computation time, using approximate, aggregated MDP models is more suitable. However, simple heuristics such as myopic policies can also provide good solutions depending on the server capacity and costs, and can be used when the computation time matters. 

The remainder of this paper is organized as follows. In Section \ref{Sec:Lit}, we present the literature review on related problems and solution methodologies. In Section \ref{Sec:Definition}, we present the problem definition and our MDP model. Then in Section \ref{Sec:Methodology}, we present our solution methodology and structural results. Section \ref{Sec:Comp} provides our computational results evaluating the performance of the proposed solution approaches. The paper is concluded in Section \ref{Sec:Conclusion} with a brief discussion on our key findings.

\section{Literature Review}
\label{Sec:Lit}

Our work is related to study of admission control, order acceptance, scheduling, and dynamic capacity allocation problems in service and manufacturing systems. Many related studies consider admission, scheduling and capacity allocation decisions simultaneously, which requires novel solution techniques due to increased problem complexity.   Here we detail the solution approaches developed in different domains. We especially give close attention to healthcare and manufacturing systems where consideration of customer classes are important. 

In service systems, the problem of admission control and scheduling generally go hand in hand. For instance, the basic patient scheduling problem consists of simultaneous admission and scheduling of patients. Generally, it is assumed that the service provider can only decide the admission of outpatients, while emergency patients must always be accepted immediately. That is, admission and scheduling decisions can only be made for one class of arrivals, which is similar to our setting. \cite{Kolesar1970} provides an early discussion of when to schedule elective patients in emergency patients' presence when the number of beds is limited. In \cite{Kolesar1970}, several queueing models and a Markovian decision model are presented. \cite{Ayvaz2010} consider the joint of admission control and dynamic allocation of hospital capacity among emergency and elective patients. Differing from \cite{Kolesar1970}, they allow to lose emergency patients at a cost. They show the monotonicity properties of the optimal policy, and propose a simple threshold heuristic policy.

The service system we consider shows similarities with a typical job shop capable of processing a wide range of products. Hence, we also review studies that consider order acceptance and capacity allocation in manufacturing with multiple job classes and uncertain arrivals. \cite{Slotnick2011} provides an overview of commonly used approaches. In practice, it might seem preferable to postpone orders until a certain time; but studies show that early processing of some jobs can be beneficial under some capacity restrictions. In an earlier study, \cite{Slotnick1996} proposed a single-period order acceptance model to maximize total profit with weighted lateness, where delayed order delivery incurs a penalty and earliest delivery contributes to profit. Since the objective does not only consider the lateness cost, the model successfully reflects the manufacturer's trade-offs. \cite{Lewis2002} further extend the former study to a multi-period case where the job selection affects future orders. \cite{Duenyas1995} formulates the problem of sequencing a single product in a manufacturing system with multiple classes of customers as a semi-Markov decision process and examines the due date setting policy. Following this idea, \cite{Carr2000} address a combined problem of admission control and scheduling in a production system with two product classes. The first class of products is made-to-stock, and the producer has obligations for their production. However, the second class of products is made-to-order, and the firm has the option of admission or rejection of a particular order. In this aspect, this study is similar to ours.

A closely related capacity allocation problem is studied by \cite{Cavdar2020} in the context of curbside pickup systems. They consider a service system with only one type of customers and no rejections. They show that for systems with short customer order horizons, a class of threshold policies are optimal for the capacity allocation; and, they develop heuristics for general systems based on threshold structures. Our problem setting is more complex compared to \cite{Cavdar2020} since we introduce different customer classes and also consider admission decisions. Our study also differs in terms of solution approaches.

MDPs have been widely used to model admission control and dynamic capacity allocation problems. Many studies that focused on patient admission and scheduling problems developed MDP models. \cite{Green2006} focus on a hospital's magnetic resonance imaging (MRI) center, providing service for different patient groups. Their analysis shapes around two interrelated tasks: designing the outpatient appointment schedule and dynamic priority rules for admitting the patient into service. They formulate the capacity allocation and the dynamic capacity allocation problem as a finite-horizon MDP. \cite{Patrick2008} explore an advanced patient scheduling problem also inspired by diagnostic facilities. In their setting, it is assumed that arriving patients can be rejected, scheduled in available slots on the same day, diverted, or postponed to the next day. They model the scheduling problem as an infinite-horizon discounted MDP and use linear-programming (LP)-based approximate dynamic programming (ADP) to solve it. \cite{Gocgun2014} study the chemotherapy scheduling problem where patients of different types dynamically arrive over time and have specific target dates along with tolerances. They model the problem as an MDP and use LP-based ADP to obtain approximate solutions. In these last two studies, there is an infinite time horizon but a finite rolling booking horizon, whereas \cite{Green2006} consider a particular time frame, e.g., a day. Our study has similarities with \cite{Patrick2008} and \cite{Gocgun2014}, both papers also consider patient priorities. In addition, \cite{Gocgun2014} integrate admission control with patient preferences which is again similar to our modeling framework. 

We obtain the optimal admission and capacity allocation policies by solving the corresponding LP for our MDP model. Solving an MDP becomes particularly challenging for models with large state and action spaces (\cite{Bhulai}), as in our case. For such problems, action elimination methods can be used to reduce the cardinality of the policy space. \cite{Lasserre1994} yields sufficient conditions that permit the detection of optimal and non-optimal actions using LP formulation. However, this method does not permanently eliminate non-optimal actions from the action-space; and it needs internal tests in various stages of the solution algorithm. Our approach differs from (\cite{Lasserre1994}) since we propose an offline action elimination method based on problem properties. 

In summary, the problem we investigate and the proposed solutions are different than the previous studies in the following major ways: $(i)$ We provide an MDP formulation for the joint admission control and capacity allocation problem in service systems with limited capacity, as well as earliness and rejection costs. In addition, our model uses information available on the arrival process to inform capacity management decisions. $(ii)$ We provide an offline action elimination method that is easy to apply since it is based on system parameters and independent of the random arrival process. Thus, it can be easily integrated with different solution approaches. $(iii)$ We develop a problem-specific aggregation method that finds approximate-optimal policies, and provide a sensitivity analysis on both the resulting optimality gap and computational performance.
\section{Problem Definition and Formulation}
\label{Sec:Definition}

We consider a service system with constant service capacity $M$ and a fixed customer order horizon length. We assume that customers arrive at each discrete time period with a service request for the current or a future period. That is, a customer arriving at time $t$ can make a request to be served in period $t+j$, where $0\leq j < K$, where $K$ represents the length of the customer order horizon. We consider two types of customers: high-priority and low-priority. We use index $i=1$ for the former and $i=2$ for the latter. High-priority customer requests must be admitted into the system. However, the service provider can reject low-priority requests at the time of arrival at a cost $c_r$. Once a low-priority request is admitted into the system, it has to be served. Customers who enter the system do not abandon. Considering the practical implications (and also the cost) of entering the system early from the customer perspective, we assume that low-priority customers are considered to be high-priority if they already waited one period to receive service. That is, their priority level is shifted from low to high, and they are treated as high-priority customers thereafter.

Each service request has to be served by its requested delivery time, i.e., while early service is allowed, late service is not. If the number of requests that are due in the current period exceeds the available capacity, we assume that these requests can be served by working in overtime or outsourced at an additional cost $c_o$ per request. 
In case of early service, the server incurs early service cost per period per request. This can be considered as the holding/handling cost of the order until pickup or customer dissatisfaction due to receiving service prior to the requested time. In case of latter, it is reasonable to expect the early service cost to be higher for high-priority customers. Thus, we distinguish the early service costs by priority class and let $c_i^e$ be the early service cost per period for priority class $i$. 
Our goal is to identify simultaneous admission and capacity allocation decisions to serve both customer classes in order to minimize long-run average operation costs.


We model this system as a discrete-time infinite-horizon MDP under the following set of simplifying assumptions. The service duration of each job is assumed to be one period, and the delivery is instantaneous; hence, we use service and delivery interchangeably throughout the paper. Decision epochs are just after customer arrivals at each discrete time period, which means the process is observed by the controller at these discrete time points when the arrivals are realized. Due to these assumptions, at the beginning of each time period, the service capacity is $M$. Therefore, we do not need to keep track of service capacity in the system state. Arriving requests are first placed in a queue named \textit{acceptance queue}. If a job is admitted into the system, it is then placed in a \textit{processing queue}. Processing queue lengths show the number of jobs of each priority waiting to be served associated with their service request times. At each decision epoch, the controller observes the acceptance and processing queues, and simultaneously decides how many low-priority jobs to admit from the acceptance queue and which jobs from processing queues to serve. Therefore, the state of the system can be described by the lengths of the acceptance and processing queues at time period $t$. We let $x^{t} = x $ be the vector denoting the length of processing queues at time $t$ immediately prior to admission and service decisions, with $x=\{x_{i j}\ |\ i=1,2\ \text{and} \ 0 \ \leq \ j \ < K\}$, where $x_{i j}$ is the number of  class $i$ jobs in the $j^{t h}$ processing queue that requested to be served $j$ periods later. Note that these are the jobs that are already accepted to the system and have to be served. Let $a^{t}=a$ with $a=\{a_{i j}\ |\ i=1,2\ \text{and} \ 0 \leq j < K\} $ represent the arrivals at time period $t$, where $a_{i j}$ is the number of class $i$ arrivals requesting to be served $j$ periods later. We represent the system state $s^{t}$ at time $t$ by the pair $(x^{t}, a^{t})$, and denote the set of all possible states by $\mathcal{S}$.

Customer arrivals follow a discrete probability distribution. The overall arrival rate for the total number of job requests that can be observed in one period is denoted by $\lambda$. The probability of observing $a$ arrivals of class $i$ customers requesting service for $t+j$ is computed with respect a to general arrival distribution with arrival rate $\lambda_{i j}$, $i=1,2$ and $0\le j < K$. 
The arrivals of each priority class $i$ and for each time period $t+j$ are independent of the number of jobs waiting to be served in the processing queues for all $i$ and $j$. We truncate arrival distribution to ensure that the state space is finite and avoid the curse of dimensionality. In particular, we assume that in each period, at most $A$ arrivals of class $i$ customers requesting to be served at period $t+j$ can be observed. Therefore, at most $2K A$ arrivals in total can be observed in a period.

In our MDP formulation, decision epochs are immediately after the arrivals realize. At each decision epoch, the controller simultaneously decides which job requests to reject among low-priority arrivals and which jobs to serve from the processing queues considering the server capacity, potential future arrivals and associated costs. We denote our action by $d = (r, y)$ where $r= \{r_{j} \ |\ 0\leq j < K \}$ represents the number of rejected low-priority job requests for time period $t + j$, and $y = \{y_{i j}\ |\ i=1,2\ \text{and} \ 0\le j < K \} $ where $y_{i j}$ is the number of class $i$ customers requested to be served at $t+j$ but served at time $t$. Notice that, $r_{j}\le a_{2j}$ holds for all $j$, since only newly arriving job requests can be rejected. Similarly, we have $y_{1j} \le x_{1j} + a_{1j}$ for high-priority jobs and $y_{2j} \le x_{2j} + a_{2j} - r_{j} $ for low-priority jobs for all $j$. Therefore, the action space is limited by the state space $\mathcal{S}$. We denote the set of all possible actions by $\mathcal{D}$ and denote the set of available actions in state $s$ by $\mathcal{D}_s$ . We depict the evolution of the described system in Figure \ref{Fig:Flow} in Appendix \ref{App:Flow}.

Given a current state $s^t$, an action $d$, and a realization of arrivals $a^{t+1}$ in the next period, the next state can be determined by the following transition function $f(.)$ and the associated transition probability $P(s^{t+1}\ |s^t,d \ )$.
\begin{equation}
    f(s^t,\ d,\ a^{t+1}) = s^{t+1} = ( x^{t+1}, a^{t+1} )
\end{equation}

where $x^{t+1}=\{x^{t+1}_{i j}\ |\ i=1,2\ \text{and} \ 0 \ \leq \ j \ < K\}$ is described with the following transitions
\begin{align*}
& \ x_{1 0}^{t+1}= ((x_{1 1}^t + a_{1 1}^t)-y_{1 1} ^{t}\ )+((x_{2 1 }^t+(a_{2 1}^t-r_{1}^t\ ))-y_{2 1}^t\ ), \\
&\ x_{1 j}^{t+1}=(x_{1 j+1 }^t + a_{1 j+1}^t)-y_{1 j+1}^t,\ \text{for} \ 1\leq j < K,  \\
&\ x_{2 0}^{t+1} = 0,  \\
&\ x_{2 j}^{t+1}=(x_{2 j+1}^t + (a_{2 j+1}^t-r_{j+1}^t\ ))-y_{2 j+1}^t,\ \text{for} \ 1\leq j < K.  
\end{align*}
where $x_{1 0}^{t+1}$ and $x_{2 0}^{t+1}$ are the number of high-priority and low priority jobs that are due in the current period, respectively. Similarly, $x_{1 j}^{t+1}$ and $x_{2 j}^{t+1}$ denote the number of high and  low priority jobs that are due in period $t+j$. 

As the transition function indicates, whenever an action is taken, the remaining unserved jobs in the system are placed in the processing queues to be served later. The associated transition probabilities are computed as follows:

\begin{equation}
  P(s^{t+1}\ |s^t, d \ )=
    \begin{cases}
    \ & \displaystyle \prod_{\shortstack{$\scriptstyle i=1,2$ \\ $\scriptstyle j\in[1,K-1]$}} p_{i j}(a_{i j}^{t+1}\ ),\  \text{if} \ s^{t+1}=f(s^t, d, a^{t+1}) \\
    \ & \displaystyle 0,\ \text{otherwise} 
    \end{cases}  
\end{equation}

where $p_{i j}(a_{i j}^{t+1}\ )$ is the probability of observing $a_{i j}^{t+1}$ customers of class $i$ with a service request that is due in period $t+j$ at the very beginning of the next period. We provide the following example to demonstrate the transitions between states.

\textbf{Example:} Consider a system with $K=3$ that is in state $s^{t}$ at time $t$, where:
\begin{equation*}
    s^{t}=(x^{t}, a^{t})=[(\underbrace{\underbrace{1,2,0}_{high},\underbrace{0,2,0}_{low}}_{processing\: queues}), (\underbrace{\underbrace{1,0,0}_{high},\underbrace{2,1,1}_{low}}_{arrival \: queue})]
\end{equation*}

In this state, the queue contains one high-priority customer who is to be served in the current period, i.e., $x_{1 0} = 1$, two high-priority and two low-priority job requests for time period $t+1$, i.e., $x_{1 1} = 2$ and $x_{2 1} = 2$, one arriving high-priority job request to be fulfilled at time period $t$, i.e., $a_{1 0} = 1$, and four low-priority arrivals in total; two requesting to be served in the current period, one requesting to be served in period $t+1$ and one requesting to be served in period $t+2$, i.e., $a_{2 0} = 2$, $a_{2 1} = 1$, $a_{2 2} = 1$. 
Assume that the controller chooses the action $d =(r,y)=[(1,0,0),(2,0,0,1,0,0)]$. In the next period, we must have $x^{t+1}=(5,0,0,0,1,0)$ before the realization of arrivals at time period $t+1$. Thus, the probability of transitioning into state $s^{t+1}=[(5,0,0,0,1,0), (0,0,1,0,2,0)]$ is
the probability of observing arrivals of one high-priority job that is due in period $t+2$ and two low-priority jobs that are due in period $t+1$. Therefore, $P(s^{t+1} | s^t, d)=P(a^{t+1} = (0,0,1,0,2,0))=P(a_{10}=2)P(a_{11}=0)P(a_{12}=0)P(a_{20}=1)P(a_{21}=0)P(a_{22}=0)$. On the other hand, the probability of transitioning into the state $[(5,0,0,0,0,0), (0,0,1,0,2,0)]$ is zero because we must have $x^{t+1}=(5,0,0,0,1,0)$ under the chosen action.

The total cost incurred in each period is the sum of the rejection, overtime, and the early service costs. For a given period, the immediate cost that is incurred under the action $d$ at state $s^t$ is denoted by $c(s^t,d)$ and is computed as follows:

\begin{equation}
   c(s^t, d)=c_o \bigg{(}M-\sum_{i=1}^{2}{\sum_{j=0}^{K-1}{y_{i j}}}\bigg{)}^+ + c_r\sum_{j=0}^{K-1} {r_{j}} + \sum_{j=1}^{K-1}{j (c_e^1 y_{1j} + c_e^2 y_{2j})} 
\end{equation}

where $(\cdot)^+$ denotes $\max(\cdot,0)$.

We aim to determine a policy that minimizes the long-run average cost. The average reward optimality criterion ultimately depends on the limiting structure of the underlying stochastic process \cite{Puterman}. Therefore, we first classify the problem. A simple analysis concludes that the resulting MDP is unichain. Regardless of the policy chosen and the initial state, we can reach the state where there are no customers in the processing queues and the acceptance queue with some positive probability since the probability of observing zero customer arrivals is positive regardless of state. Therefore, the state with zero customers is recurrent. Since recurrence is a class property, all communicating states that are accessible from zero customer state create a recurrent class.

For unichain models, all stationary policies provide constant gain \cite{Puterman}, a constant cost in our case. In addition, for unichain models with finite state and action spaces and finite immediate rewards, the existence of an optimal stationary deterministic policy is guaranteed. Therefore, a single optimality equation is enough to characterize optimal policies and their average costs. Accordingly, the optimal policy satisfies the following optimality equation \cite{Puterman}:
\begin{equation}\label{Eq::Optimality}
    0=\min_{d \in \mathcal{D}_s}\bigg{\{}c\left(s,d\right)-g+\sum_{s^\prime\in \mathcal{S}}P(s^\prime|s,d)h(s^\prime)-h(s) \bigg{\}}, \quad \forall s \in \mathcal{S}
\end{equation}

where $g$ denotes the long-run average cost and $h(\cdot)$ is the bias.


The dimension of the resulting MDP rapidly increases with the number of time periods in the customer order horizon $K$ and the maximum number of arrivals $A$ that can be observed in a period for each customer class. Solving MDPs with large and multidimensional state and action spaces is computationally challenging, especially if one is using standard iterative methods, such as value iteration. Therefore, we direct our efforts to developing efficient solution methods. In this paper, we propose a linear programming (LP)-based solution method. To increase the computational performance, we partially characterize the suboptimal solutions. In addition, we propose a problem-specific approximation method through a state space reduction technique.
\section{Solution Methodology}
\label{Sec:Methodology}
In this section, we present several methods for solving the joint admission and capacity allocation problem in a two-class service system. We first present the LP formulation of the described problem as follows:
\begin{maxi!}|s|[2]
    {}                               
    { g \label{peq:eq1}}   
    {\label{Primal}}             
    {}                               
    \addConstraint{g+h(s)-\sum_{s^\prime\in \mathcal{S}}P(s^\prime|s,d)h(s^\prime)}{ \leq c(s,d), \forall d \in \mathcal{D}_s, \forall s \in \mathcal{S} \label{peq:con1}}    
\end{maxi!}

We work with the dual of this formulation where the decision variables $x(s,d)$  are defined for each state and each action that are allowable in that state. 
\begin{mini!}|s|[2]
    {}                               
    { \sum_{s \in \mathcal{S}} \sum_{d \in \mathcal{D}_s} c(s,d)x(s,d) \label{deq:eq1}}   
    {\label{Dual}}             
    {}                               
    \addConstraint{ \sum_{d \in \mathcal{D}_{s^\prime}}x(s^\prime, d) - \sum_{s \in \mathcal{S}} \sum_{d \in \mathcal{D}_s} P(s^\prime|s,d)x(s,d)}{ = 0, \forall s^\prime \in \mathcal{S} \label{deq:con1}}    
    \addConstraint{\sum_{s \in \mathcal{S}} \sum_{d \in \mathcal{D}_s}x(s,d)}{ = 1 \label{deq:con2}}  
    \addConstraint{x(s,d)}{\ge 0, \forall d \in \mathcal{D}_s, \forall s \in \mathcal{S} \label{deq:con3}}
\end{mini!}

In the dual LP (DLP) model, each feasible solution corresponds to a stationary policy and the decision variable $x(s,d)$ can be interpreted as the stationary probability of being in state $s$ and choosing action $d$ under the policy. Moreover, each basic feasible solution to the DLP corresponds to a stationary deterministic policy $\pi$ such that $\pi_x(s)=d$ if $x(s,d)>0$. Note that unichain models have possibly an empty set of transient states, then $\pi_x(s)$ can be chosen arbitrarily for transient state $s$ with $x(s,d)=0$ for all $d \in \mathcal{D}_s$.

Our MDP model grows exponentially in the length of customer order horizon and the maximum number of arrivals that can be observed in a period. As a result, using the DLP model directly creates significant computational limitations to solve the problem, even though the complexity of LP approach is known to be polynomial for MDPs  \cite{Papadimitriou1987}. Therefore, to increase the computational efficiency of the DLP model, we propose a set of action and state space reduction methods, called \textit{action elimination} and \textit{state aggregation}, respectively. To reduce the size of the action space, we partially characterize the suboptimal actions using some structural results and successfully restrict the policy space, i.e., the search space of the corresponding DLP model. Hence, the presented action elimination procedure ensures a significant improvement in the performance of the DLP model in terms of solution time and memory requirement. We also propose a state aggregation method, and thus an aggregate MDP model. We then integrate the action elimination and the state aggregation procedure for further improvements.

\subsection{Action Elimination}

In this section, we partially characterize the suboptimal policies by determining some necessary conditions that a set of feasible actions for a state $s^t \in \mathcal{S}$ should satisfy for optimality. To avoid cases with trivial solutions, we make two assumptions on the cost parameters. Assumption \ref{Assumption::1} reflects the operational costs of early service option and Assumption \ref{Assumption::2} ensures that the rejection decision remains relevant in the problem.

\begin{assumption}
\label{Assumption::1}
	$c_{2}^e < c_{1}^e$: serving a high-priority job request early is costlier than serving a low-priority job request early. 
\end{assumption}
\begin{assumption}
\label{Assumption::2}
	$c_{r} < c_{o}$: rejection of a low-priority job request incurs less cost than its fulfillment in overtime. 
\end{assumption}
\begin{proposition}
\label{Proposition::1}
    For arbitrary cost parameters, early service is suboptimal for a job request in the processing queue unless there is remaining capacity after serving job requests that are due in the current period.
\end{proposition}

The proof of the proposition is in Appendix \ref{Ap:Pros}. Proposition \ref{Proposition::1} implies that any policy  that uses an action that serves a job early in a state $s$ using over time is suboptimal. Thus, such actions can be eliminated from the set of allowable actions $\mathcal{D}_s$. In the following part, we present our structural analysis. Our initial results are derived for systems with $K=2$ and summarized in Proposition \ref{Proposition::2}. 

\begin{lemma}
When $K=2$, if an early service decision should be taken, the service provider should start with the low-priority jobs first.
\end{lemma}

This lemma follows Assumption \ref{Assumption::2}. When the customer order horizon has two periods, if a job that is due in period $t + 1$ is served early, we incur $c_{1}^e$ or $c_{2}^e$ for a high-priority and a low-priority job, respectively. Since $c_{2}^e < c_{1}^e$, whenever early service is reasonable, the controller should start from low-priority jobs.

\begin{proposition}
\label{Proposition::2}
	When $K=2$, all low-priority arrivals that would need to be served in overtime in the current period or that would cause overtime in the next period if they are admitted should be rejected. Hence, an action $d=(r,\ y)$ in a state $s\in \mathcal{S}$ is suboptimal if it does not satisfy the following conditions:
\begin{itemize}
 \item When $c_{2}^e < c_{1}^e < c_{r} < c_{o}$: 
 \vspace*{-5mm}
	\begin{align*} 
    & (i)\ min\{ [(x_{10}+a_{10}+a_{20})-M]^+,a_{20} \}  \le r_{20}\le a_{20} \\
    & (ii)\ min\{[[(a_{11}+a_{21}\ )-M]^+-y_{21} ]^+,a_{21}  \} \le r_{21} \le     [a_{21}-[M-(x_{10}+a_{10}+a_{20}-r_{20}\ )]^+ ]^+ \\
    & (iii)\ y_{11} =
       \begin{cases}
     0 \ \text{if} \ y_{21}<min\{[M-(x_{10}+a_{10}+a_{20}-r_{20})]^+,(a_{21}-r_{21})\} \\
    \le min\{[M-(x_{10}+a_{10}+a_{20}-r_{20}+y_{21})]^+,a_{11} \}, \ \text{otherwise}
    \end{cases} \\
    & (iv)\ y_{21} \le min\{[M-(x_{10}+a_{10}+a_{20}-r_{20})]^+,(a_{21}-r_{21} )\}
    \end{align*}
\item When $ c_{r} < c_{2}^e < c_{1}^e < c_{o}$: 
 \vspace*{-5mm}
\begin{align*}
   & (i)\ min\{ [(x_{10}+a_{10}+a_{20})-M]^+,a_{20} \}  \le r_{20}\le a_{20}\\
   & (ii)\ min\{[(a_{21}+a_{21})-M]^+,a_{21} \} \le r_{21} \le a_{21} \\
   & (iii)\ y_{11} \le min\{[M-(x_{10}+a_{10}+a_{20}-r_{20}+y_{21})]^+,a_{11} \} \\
   & (iv)\ y_{21} = 
   \begin{cases}
    0,\ \ \text{if}\ r_{21}>0 \\
    \le min\{[M-(x_{10}+a_{10}+a_{20}-r_{20})]^+, (a_{21}-r_{21} )\} , \ \text{otherwise} \
   \end{cases}
\end{align*}
\item When $ c_{r} < c_{o} < c_{2}^e < c_{1}^e $: 
 \vspace*{-5mm}
    \begin{align*}
   & (i)\ min\{ [(x_{10}+a_{10}+a_{20})-M]^+,a_{20} \}  \le r_{20}\le a_{20}\\
   & (ii)\ min\{[(a_{21}+a_{21})-M]^+,a_{21} \} \le r_{21} \le a_{21} \\
   & (iii)\ y_{11} = 0 \\
   & (iv)\ y_{21} = 0
    \end{align*}
\end{itemize}
 \vspace*{-5mm}

with $y_{11}+y_{21}\le [M-(x_{10}+a_{10}+a_{20}-r_{20})]^+$.
\end{proposition}

The proof of the proposition is presented in Appendix \ref{Ap:Pros}. In the proposition, we do not provide an expression for the actions that must be taken for the jobs that are due in the current period, i.e. $y_{1 0}$ and $y_{2 0}$. These actions follow from our problem definition; all jobs that are due the current period, should be fulfilled.  

These results identify a set of suboptimal policies based only on system parameters such as server capacity $M$ and costs without relying on the arrival probabilities. Therefore, these results enable us to eliminate a significant number of the suboptimal policies quickly from the policy space for systems with $K=2$. To address more general systems, we extend some of these results to larger customer order horizons. 

\begin{definition}
\label{Def::t-boundary}
	A state $s \in \mathcal{S}$ is a $j$-boundary state, if there are at least $M$ jobs that are due $t+j$, $j\in \{0,\ldots,K-1\}$ in state $s$, including the jobs in the processing queue and the acceptance queue prior to any decision.
\end{definition}

\begin{proposition}
\label{Proposition::3}
An action $d=(r,\ y)$ is suboptimal in a state $s\in \mathcal{S}$ that is a $j$-boundary state for all $j\in \{0, 1, \dots, K-1\}$ if it does not reject at least as many jobs as the number of jobs that will cause overtime for the next $K$ periods including the current period, and it serves at least one job early regardless of its priority.
\end{proposition}



The proof is in Appendix \ref{Ap:Pros}. Proposition \ref{Proposition::3} is an immediate result on suboptimal actions for arbitrary $K$ and used as a first step in action elimination without any further considerations. The following proposition is the extension of Proposition \ref{Proposition::2} for longer customer order horizon. Note that for longer customer order horizons, the relations between costs, e.g., $c_r/c_{2}^e$, become more influential. 

\begin{proposition}
\label{Proposition::4}
An action $d=(r,\ y)$ is suboptimal for state $s^t \in \mathcal{S}$ if it does not satisfy the following conditions. 
    \begin{itemize}
        \item  When $c_{2}^e < c_{1}^e < c_{r} < c_{o}$:
         \vspace*{-5mm}
            \begin{align*}
   & (i)\ min\{ [(x_{10}+a_{10}+a_{20})-M]^+,a_{20} \}  \le r_{20}\le a_{20}\\
   &
   \begin{aligned}
    & (ii)\  min\{[overtime_{j}-y_{2 j}]^+, a_{2 j} \}  \le r_{2 j} \le a_{2 j} \ \mbox{for} \  j=1,\ldots,K-1 \ \mbox{if} \ j c_{2}^e \le c_r \\
    & (iii)\ min\{overtime_{j},a_{2 j} \}  \le r_{2 j} \le a_{2 j} \ for \  j=1,\ldots,K-1 \ \mbox{if} \ j c_{2}^e > c_r
   \end{aligned} \\
   & (iv)\ y_{1 j} \le min\{remain_j ,x_{1 j} + a_{1j} \}, \ \mbox{for} \  j=1,\ldots,K-1 \\
   & (v)\ y_{2 j} \le min\{remain_j,(x_{2j}+a_{2j}-r_{2j} )\} \ \mbox{for} \ j=1,\ldots,K-1
\end{align*}
\item When $c_{r} < c_{2}^e < c_{1}^e < c_{o}$:
         \vspace*{-5mm}
\begin{align*}
   & (i)\ \ min\{ [(x_{10}+a_{10}+a_{20} )-M]^+,a_{20} \}  \le r_{20}\le a_{20} \\
   & (ii)\ min\{overtime_j,a_{2j} \}  \le r_{2j} \le a_{2j} \ \mbox{for} \  j=1,\ldots,K-1 \\
   & (iii)\ y_{1j} \le min\{remain_j , (x_{1 j} + a_{1 j}) \}, \ \mbox{for} \  j=1,\ldots,K-1 \\
   & (iv)\ y_{2j} \le min\{remain_j,(x_{2j}+a_{2j}-r_{2j} )\} \ \mbox{for} \ j=1,\ldots,K-1
\end{align*}
\item When $ c_{r} < c_{o} < c_{2}^e < c_{1}^e $:         \vspace*{-5mm}
\begin{align*}
   & (i)\ min\{ [(x_{10}+a_{10}+a_{20} )-M]^+,a_{20} \}  \le r_{20}\le a_{20} \\
   & (ii)\ min\{[(x_{1j}+x_{2j}+a_{1j}+a_{2j})-M]^+,a_{2j} \}  \le r_{2j} \le a_{2j} \ for \  j=1,\ldots,K-1 \\
   & (iii)\ y_{1j}=0 \ for \ j=1,\ldots,K-1 \\
   & (iv)\ y_{2j}=0 \ for \ j=1,\ldots,K-1
\end{align*}
    \end{itemize}
            \vspace*{-5mm}
where $overtime_{j}= [[(x_{1j}+x_{2j}+a_{1j}+a_{2j})-M]^+$, $remain_j=[M-(y_{10} + y_{20} +\sum_{j=1}^{j-1}{y_{2 j}} + \sum_{k=0}^{j-1}{y_{1 k}} ) ]^+$ with $\sum_{j=1}^{K-1}{y_{1j} + y_{2j}} \le [M-(x_{10}+a_{10}+a_{20}-r_{20})]^+ $.

\end{proposition}

The proof of the proposition is in Appendix \ref{Ap:Pros}. 

When constructing the DLP model associated with our MDP model, we first eliminate the suboptimal actions from the allowable action set for each state based on the results developed in this section. This procedure decreases the size of the policy space significantly and strengthens the DLP model. In the following section, we propose a state aggregation method to decrease the computational requirements further.

\subsection{Aggregate MDP Model}

In this section, we introduce methods that can be used to aggregate the state space to find fast heuristic solutions. The use of aggregation techniques is a common approximation approach. The key principle of this approach is to cluster the original large state space into aggregate subgroups, which are treated as newly created states in the aggregate MDP. Aggregation decreases the sizes of the state space and the transition probability matrix. Using this approach, the original MDP can be approximated by a smaller aggregate MDP. Since the aggregate model preserves the Markov property, any solution method for MDP's, such as LP, is still valid. Solving the aggregate MDP model, we can find approximate solutions for the original MDP that can be evaluated on the original space, where the level of aggregation determines the quality of solutions.

To aggregate the state space, we follow a method that is commonly used in the MDP literature, more specifically for MDP applications in reinforcement learning \cite{sutton2018reinforcement}. Our aggregation is based on clustering the states that are considered to be ``equivalent" according to a certain criterion into \textit{meta-states}. Equivalent states are assumed to behave approximately the same under all policies. Based on the equivalence definition by \cite{Givan2003}, we define $\gamma$-equivalence below.


\begin{definition}
\label{Def::2}
   Given an MDP, described with tuple $< \mathcal{S}, \mathcal{D}, P, c >$, two states $s, s^{\prime} \in \mathcal{S}$ with the set of actions $\mathcal{D}_s, \mathcal{D}_{s^{\prime}} \subset \mathcal{D}$ respectively, are said to be $\gamma$-equivalent if for all $d \in \mathcal{D}_s$ there is a $d^{\prime} \in \mathcal{D}_{s^{\prime}}$ with $\gamma$  such that
    \begin{equation}
        \abs{c(s, d)-c(s^{\prime}, d^{\prime})} \leq \gamma\quad \text{and} \quad  P(s^{\prime\prime} | s, d)= P(s^{\prime \prime} | s^{\prime}, d^{\prime}), \ \forall s^{\prime\prime} \in \mathcal{S}
    \end{equation} 
\end{definition} 
In Definition \ref{Def::2}, the value of $\gamma$ affects the level of aggregation. Smaller values of $\gamma$ decrease the level of aggregation and potentially result in better solution quality, but increase the computation time of the aggregation. Larger values, on the other hand, speed up the aggregation procedure, but decrease the solution quality due to higher level of aggregation. Therefore, this trade-off should be considered when choosing the value of $\gamma$. In addition, while a static value of $\gamma$ can be used, it can also be updated dynamically during the aggregation process. 

Before aggregating the states according to Definition \ref{Def::2}, we enforce an additional condition to increase the similarity of states that are clustered together. We aggregate two states $s=(x,a)$ and $s'=(x', a')$ according to Definition \ref{Def::2} as long as $x_j + a_j = x^{\prime}_j + a^{\prime}_j$ for all $j$, $0\leq j <K$. Under this condition, states $s$ and $s'$ have the same number of jobs; however, they are not the same state unless they have the same acceptance queues $a$ and $a'$. We call this aggregation method \emph{total-job aggregation}. The example below demonstrates a case where two states are clustered together.  

\textbf{Example:} Consider a state $s=(x,a)=[(1,2,0,0,2,0), (1,0,0,2,1,1)]$ for $K=3$ and $A=2$. Let the service capacity be $M=3$. Then, $s$ is a $0$-boundary state, since $x_{10}+a_{10}+ x_{20} + a_{20} > 3$. Therefore, an optimal action in $s$ can not perform early service. Now consider a second state $s^{\prime}=(x^{\prime},a^{\prime})=[(2,2,0,0,1,0), (0,0,0,2,2,1)]$ under the same problem settings. $s^{\prime}$ is also a $0$-boundary state, and $x_j + a_j = x^{\prime}_j + a^{\prime}_j$ for all $j$, $0\leq j <3$. The action $d=(r,y)=[(1,0,0), (3,0,0,0,0,0)]$ is allowed in both $s$ and $s'$, leads to the same state in the next period, and incurs the same cost for both. Therefore, states $s$ and $s'$ are aggregated into the same meta-state according to our total-job aggregation method.

The aggregation method clusters the original state space into meta-states $ \mathbf{S}=\{S_1, \ldots, S_N \}$, where $N$ is the number of meta-states. Note that, for each $i \in \{1, \ldots, N\}$, the set of actions available in state $S_i$ in the aggregate MDP, $\mathcal{D}_{S_i}$ is the intersection of the corresponding sets of allowable actions in the original MDP, i.e., $\mathcal{D}_{S_i}=\cap_{s \in S_i} \mathcal{D}_s$. We define the one-step transition probabilities and costs for the newly defined aggregate states $S_i, S_j \subseteq \mathbf{S}$ and $d \in \mathcal{D}_{S_i}$ as follows:
\begin{equation}
\label{Aggregate::Probability}
     P(S_j | S_i, d) = \sum_{s \in S_i} \frac{1}{|S_i|} \sum_{s'\in S_j}P(s'| s, d)  \quad \forall S_i, S_j \in \mathbf{S} \quad \forall d \in \mathcal{D}_{S_i}, 
\end{equation} 
\begin{equation}
\label{Aggregate::Cost}
    c(S_i, d) = \sum_{s \in S_i} \frac{1}{|S_i|} c(s,d) \quad \forall d \in \mathcal{D}_{S_i}. 
\end{equation} 

\begin{proposition}
\label{Pro:AggMarkov}
Transition probabilities for the aggregate state space, described in equation (\ref{Aggregate::Probability}), constitute a stochastic matrix.
\end{proposition}

The proof of the proposition is presented in Appendix \ref{Ap:Pros}. Proposition \ref{Pro:AggMarkov} shows that the aggregate model maintains a Markov decision process structure. The LP model can be implemented on the aggregate MDP to find heuristic solutions in shorter amount of computation time. The solution to the aggregate model is calculated using the average costs described in equation (\ref{Aggregate::Cost}). Once we determine the solution to the aggregate model, we can calculate the exact value of this policy on the original model.
Let $\bar{\pi}$ represent the optimal policy for the aggregate model. $\bar{\pi}$ dictates a policy $\delta$, for the original MDP, described as follows:
\begin{equation}
\label{Aggregation::Evaluation}
\delta(s)=\bar{\pi}(S_{i}) \text { if } s \in S_{i}, \\ \quad \forall s \in S \text{ and }\forall S_i \in \mathbf{S}.
\end{equation} 
So far, we presented our formulation and methods to increase the efficiency of the computations. In the following section, we present our computational results.

\section{Computational Results}
\label{Sec:Comp}
We present and discuss an extensive set of computational results to evaluate the performance of solution methods we developed and to investigate the impact of problem parameters on the optimal policy. Our main solution method uses the structural results discussed in Section 4.1 to identify suboptimal policies. We also consider an aggregate MDP model to find approximate solutions in shorter computation time. Therefore, we present our results in three parts. First, we compare the policies obtained by the dual linear programming (DLP) model, and our reduced linear programming (RLP) model, which is enhanced with action elimination. Since, both DLP and RLP models find optimal policies, we discuss the contribution of action elimination method in terms of computation time and computational requirements. Next, we compare the policies obtained by our aggregate MDP model to the optimal policies and the policies obtained by a benchmark aggregate model. Finally, we compare our aggregate MDP model with two benchmark heuristics. 

In our computational study, we experiment on different customer order patterns and priority distributions. We refer to customers' order pattern as the \textit{load} of the system. We use the term \textit{segmentation} to describe the distribution of arriving customers' priorities. Load and segmentation are integrated into our formulation through the probabilities $v_j$ and $q_i$, respectively. We formally define three different segmentation patterns: 
\begin{enumerate}
	\item \emph{Equal-segmented systems} (ES): Customers arriving at time $t$ are equally likely to have high-priority or low-priority, $q_i=\frac {1}{2}, \ i=1,2$.
	\item \emph{High-priority systems} (HS): Customers arriving at time $t$ are more likely to have high-priority. We choose $q_1=\frac{4}{5}$ and $q_2=\frac{1}{5}$.
	\item \emph{Low-priority systems} (LS): Customers arriving at time $t$ are more likely to have low-priority. We choose $q_1=\frac{1}{5}$ and $q_2=\frac{4}{5}$.
\end{enumerate} 

We also consider different load types in our experiments as described below:
\begin{enumerate}
	\item \emph{Equal-load systems} (EL): Customers arriving at time $t$ are equally likely to request service
	for any time period $t + j$, where $0 \leq j < K$, and $v_j=\frac {1}{K}, \ \forall j $.
	\item \emph{Front-loaded systems} (FL): Customers arriving at time $t$ are more likely to request service for time period $t + j_1$ than time period $t + j_2$ whenever $j_1 < j_2$. We calculate $v_j=\frac {(K-j)^2}{\sum_{k=1}^{K}k^2}, \ \forall j$.
	\item \emph{Back-loaded systems} (BL): Customers arriving at time $t$ are more likely to request service
	for time period $t + j_1$ than time period $t + j_2$ whenever $j_1 > j_2$. We calculate $v_j=\frac {(j+1)^2}{\sum_{k=1}^{K}k^2}, \ \forall j$.
\end{enumerate}

For the front-loaded and the back-loaded systems, we assume a quadratic change in $v_j$ values with respect to $j$. 

To model customer arrivals, we use a truncated Poisson distribution. The parameter $\lambda$ represents the overall arrival rate. We denote the arrival rate of a class $i$ customer requesting service for the time period $t+j$ by $\lambda_{i j}$. Considering the segmentation and load of the system, $\lambda_{i j}$ is calculated as $\lambda_{i j}=q_iv_j \lambda$. Based on these rates, we can calculate the truncated Poisson arrival process as follows. Note that $A$ represents the maximum number of arrivals that can be observed in a period.
\begin{equation}
	    p_{i j}(a)=\frac{\frac{e^{-\lambda_{i j}}\lambda{ i j}^a}{a!}}{\sum_{n=0}^{A}\frac{e^{-\lambda_{i j}}\lambda{i j}^n}{n!}},
	  \ \text{for} \ i=1,2 \ \text{and} \ 0\leq j < K
\end{equation}

We generate the parameters used in our experiments  by computing the Cartesian product of sets of values chosen for the problem parameters: number of servers ($M$), overall arrival rate ($\lambda$), maximum number of arriving job requests for each period ($A$), the length of the customer order horizon ($K$), the cost parameters ($c_o, c_r,  c_e^1, c_e^2$). We also consider combinations of different segmentation and the load patterns of the system. The sets of values that were used in our experiments are shown in Table \ref{Table::1}.

\begin{table}[H]
    \centering
    \caption{Parameter values used in the computational experiments.}
    \begin{tabular}{lcl}
	\hline
	$M$ & $\in$ & $\{1, 2, 5\}$ \\ 
	$A$ & $\in$ & $\{1, 2, 3, 4, 5, 6\}$ \\ 
	$\lambda$ & $=$ & $0.5A$ \\
	$K$ & $\in$ & $\{2, 3, 4\}$ \\
	$c_o$ & $=$ & $200$ \\
	$c_r$ & $\in$ & $\{50, 100, 150\}$ \\
	$c_e^1$ & $\in$ & $\{100, 150, 200, 300\}$ \\
	$c_e^2$ & $\in$ & $\{50, 150, 250\}$ \\
	\hline
    \end{tabular}
\label{Table::1}
\end{table}

Some of the parameter combinations result in very large state spaces and require computational power that exceeds the available resources. Therefore, such combinations are omitted in the following subsections. We implement all algorithms in C++ programming language and use IBM ILOCPLEX 12.10 as the solver. All computations are performed on Intel(R) Core(TM) i7 Processor 3.10 GHz with 16.0 GB RAM. For all computations, we impose a two-hour time limit.

\subsection{Comparison of the Reduced Linear Program to the Dual Linear Program}

In the first set of computational experiments, we compare the optimal solutions computed by the DLP model and the RLP model. Recall that the number of states determines the number of rows, and the number of actions determines the number of columns in the DLP formulation. Due to the elimination of a set of suboptimal actions, the RLP model has fewer columns than the DLP model. For both models, we are restricted to instances where the memory requirements are within our computational limits. In Table \ref{parameters}, we present the system parameters and the resulting number of states for the instances that can be solved to optimality. The number of states of each MDP model is determined by the length of the customer order horizon, $K$ and the maximum number of arrivals in each period, $A$. The number of actions increases in system capacity $M$ in addition to $K$ and $A$.  

\begin{table}[h]
\setcellgapes{1pt}\makegapedcells
\caption{Sizes of optimally solvable MDP models.}\label{parameters}
\centering
\resizebox{0.5\textwidth}{!}{
\begin{tabular}{r|r|r|rrr}
\hline \hline
$K$                & $A$ & \makecell{Number of \\states} & \multicolumn{3}{c}{Number of actions} \\ \hline
                   &     &                  & $M=1$      & $M=2$      & $M=5$       \\ \hline
\multirow{5}{*}{2} & 1   & 48               & 118        & 145        & 215         \\
                   & 2   & 405              & 1683       & 1896       & 3646        \\
                   & 3   & 1792             & 11416      & 12210      & 21374       \\
                   & 4   & 5625             & 51175      & 53290      & 81760       \\
                   & 5   & 14256            & 175806     & 180435     & 248032      \\ \hline
\multirow{2}{*}{3} & 1   & 1280             & 4968       & 7240       & 21497       \\
                   & 2   & 59049            & 497664     & 607095     & 2169129     \\ \hline
\multirow{1}{*}{4} & 1   & 64512            & 386640     & 653112     & 3731748   \\ 
 \hline \hline
\end{tabular}}
\label{Table::2}
\end{table}


The optimal cost and computation times strongly depend on the length of the customer order horizon $K$, service capacity $M$, maximum number of arrivals per period for each customer type $A$, costs, and segmentation and load of the system. Therefore, we organize our discussion around these components. We present the first set of results for $K=2$, $c_o=200$, $c_r=200$, $c_e^1=$, $c_e^2=20$, $A\in \{2,3,4\}$ with three types of load and segmentation in Table \ref{Table::Cost-a}. We let either $M=2$ or $M=5$ to represent systems with small and moderately large capacity. The cost of the optimal policy (obtained by both DLP and RLP models) and the computation time of each model are reported. Table \ref{Table::Cost-aComptreq}, presents the computational requirements for the corresponding MDP models.

\begin{table}[h]
\caption{Computational results for $K=2$, $c_o=200$, $c_r=150$, $c_e^1=100$, $c_e^2=50$.}
\centering
\resizebox{0.9\textwidth}{!}{
\begin{tabular}{lr|rrr|rrr|rrr}
\hline \hline
                    & \multicolumn{1}{l|}{} & \multicolumn{9}{c}{$M=2$}                     \\ \hline
                    & \multicolumn{1}{l|}{} & \multicolumn{3}{c|}{Low-priority}                                                     & \multicolumn{3}{c|}{Equal-segmented}                                                & \multicolumn{3}{c}{High-priority}                                                    \\ \hline
                    & A                     & \multicolumn{1}{c}{Cost} & \multicolumn{1}{c}{$DLPt(sec)$} & \multicolumn{1}{c|}{$RLPt(sec)$} & \multicolumn{1}{c}{Cost} & \multicolumn{1}{c}{$DLPt(sec)$} & \multicolumn{1}{c|}{$RLPt(sec)$} & \multicolumn{1}{c}{Cost} & \multicolumn{1}{c}{$DLPt(sec)$} & \multicolumn{1}{c}{$RLPt(sec)$} \\ \hline
\multirow{3}{*}{EL} & 2                     & 10.83                    & 0.71                        & 0.65                         & 13.05                    & 0.70                        & 0.59                         & 13.01                    & 0.67                        & 0.61                        \\
                    & 3                     & 35.52                    & 14.89                       & 7.34                         & 40.56                    & 14.92                       & 6.41                         & 44.35                    & 14.27                       & 5.80                        \\
                    & 4                     & 74.32                    & 193.41                      & 83.55                        & 81.85                    & 214.92                      & 67.40                        & 90.77                    & 215.90                      & 79.56                       \\ \hline
\multirow{3}{*}{FL} & 2                     & 8.93                     & 0.68                        & 0.60                         & 11.89                    & 0.65                        & 0.62                         & 10.28                    & 0.67                        & 0.58                        \\
                    & 3                     & 34.98                    & 12.58                       & 6.15                         & 40.80                    & 14.31                       & 7.16                         & 41.36                    & 12.64                       & 5.42                        \\
                    & 4                     & 73.03                    & 371.02                      & 63.59                        & 82.94                    & 224.07                      & 85.39                        & 89.36                    & 174.32                      & 52.05                       \\ \hline
\multirow{3}{*}{BL} & 2                     & 7.98                     & 0.73                        & 0.63                         & 10.09                    & 0.68                        & 0.63                         & 9.17                     & 0.67                        & 0.66                        \\
                    & 3                     & 30.47                    & 13.63                       & 6.22                         & 35.49                    & 14.23                       & 8.04                         & 36.29                    & 15.43                       & 7.71                        \\
                    & 4                     & 66.13                    & 187.19                      & 69.39                        & 75.23                    & 235.25                      & 137.20                       & 81.64                    & 227.74                      & 97.78                       \\ \hline
                    & \multicolumn{1}{l|}{} & \multicolumn{9}{c}{$M=5$}                                                                                                                                                                                                                                           \\ \hline
                    & \multicolumn{1}{l|}{} & \multicolumn{3}{c|}{Low-priority}                                                     & \multicolumn{3}{c|}{Equal-segmented}                                                & \multicolumn{3}{c}{High-priority}                                                    \\ \hline
                    & A                     & \multicolumn{1}{c}{Cost} & \multicolumn{1}{c}{$DLPt(sec)$} & \multicolumn{1}{c|}{$RLPt(sec)$} & \multicolumn{1}{c}{Cost} & \multicolumn{1}{c}{$DLPt(sec)$} & \multicolumn{1}{c|}{$RLPt(sec)$} & \multicolumn{1}{c}{Cost} & \multicolumn{1}{c}{$DLPt(sec)$} & \multicolumn{1}{c}{$RLPt(sec)$} \\ \hline
\multirow{3}{*}{EL} & 2                     & 0.01                     & 1.20                        & 1.15                         & 0.03                     & 1.06                        & 1.04                         & 0.01                     & 1.20                        & 1.18                        \\
                    & 3                     & 0.45                     & 26.18                       & 21.49                        & 0.66                     & 20.72                       & 18.25                        & 0.50                     & 26.43                       & 22.51                       \\
                    & 4                     & 2.72                     & 954.72                      & 302.45                       & 3.19                     & 924.27                      & 487.45                       & 3.15                     & 798.42                      & 260.57                      \\ \hline
\multirow{3}{*}{FL} & 2                     & 0.00                     & 1.20                        & 1.16                         & 0.01                     & 1.07                        & 1.03                         & 0.00                     & 1.15                        & 1.09                        \\
                    & 3                     & 0.24                     & 27.07                       & 19.58                        & 0.45                     & 20.14                       & 17.26                        & 0.25                     & 22.83                       & 17.29                       \\
                    & 4                     & 1.91                     & 359.76                      & 231.99                       & 2.84                     & 483.81                      & 246.52                       & 2.07                     & 423.00                      & 264.17                      \\ \hline
\multirow{3}{*}{BL} & 2                     & 0.00                     & 1.21                        & 1.19                         & 0.01                     & 1.19                        & 1.17                         & 0.01                     & 1.15                        & 1.13                        \\
                    & 3                     & 0.24                     & 24.08                       & 19.84                        & 0.46                     & 30.12                       & 21.69                        & 0.28                     & 24.18                       & 20.60                       \\
                    & 4                     & 1.76                     & 535.56                      & 238.35                       & 2.44                     & 1087.46                     & 263.50                       & 1.95                     & 1066.22                     & 253.09                      \\ \hline \hline
\end{tabular}}
\label{Table::Cost-a}
\end{table}

Table \ref{Table::Cost-a} shows a sharp increase in the computation time as $A$ increases due to the larger state space. However, when we compare the computation times of RLP and DLP models, we see that our RLP model finds the optimal policy in significantly shorter amount of time. For example, when $A=4$ and $M=5$, computation time is reduced by 50\%-75\%  due to action elimination. This indicates that the RLP model is solved on a policy space which is significantly smaller than the original. 

The impact of system load on the performance of the optimal policy can be explained by the rejection decisions. Without the rejection actions, the cost of the optimal policy would have increased as the system becomes more back-loaded (\cite{Cavdar2020}). However, when the service provider rejects low-priority jobs, the effective load in the system differs from the load of the arrivals. 
We observe that equal-load systems have higher cost than the front-loaded and back-loaded systems in general. In front-loaded systems, the rejection decisions can be done with more precision, and costly overtime can be avoided. In back-loaded systems, the service provider has more flexibility through early service to decrease the risk of overtime.
  
The impact of the customer segmentation highly depends on the service capacity. When the capacity is scarce (e.g., $M=2$ for $A=4$ in Table \ref{Table::Cost-a}), rejection and overtime costs are incurred more commonly. As a result, the cost of the optimal policy increases as the volume of high-priority customers increases since these customers have to be admitted and served even regardless of available capacity. In contrast, when the system capacity is larger (e.g.,when $M=5$), the equal-segmented systems have the highest operating costs. This is because when there is enough capacity, it is possible to admit jobs into the system that could later result in overtime due to instant overcrowding. 
In high-priority segmented systems, a large proportion of arriving customers are likely to be admitted into the system. Therefore, these systems incur higher costs compared to equal segmented and low-priority segmented systems. On the other hand, in low-priority segmented systems, a large number of customers are low priority and depending on the service capacity they can to be rejected upon arrival or served early at a smaller cost.

\begin{table}[h]
\setcellgapes{2.5pt}\makegapedcells
\caption{Computational requirements for $K=2$ when $c_o=200$, $c_r=150$, $c_e^1=100$, $c_e^2=50$.}
\centering
\resizebox{0.8\textwidth}{!}{
\begin{tabular}{r|rr|rr|rr|rr}
\hline \hline
\multicolumn{1}{l|}{} & \multicolumn{4}{c|}{$M=2$}        & \multicolumn{4}{c}{$M=5$}     \\ \cline{1-9} 
\multicolumn{1}{l|}{} & \multicolumn{2}{c|}{DLP}                                                               & \multicolumn{2}{c|}{RLP}                                                               & \multicolumn{2}{c|}{DLP}                                                                & \multicolumn{2}{c}{RLP}                                                               \\ \hline
$A$                   & \makecell{Number of \\actions} & $Mem$  & \makecell{Number of \\actions} & $Mem$  & \makecell{Number of \\actions} & $Mem$   & \makecell{Number of \\actions} & $Mem$  \\ \hline
2                     & 1896                                                                          & 4.08   & 1040                                                                          & 2.33   & 3646                                                                          & 7.69    & 3430                                                                          & 7.25   \\
3                     & 12210                                                                         & 74.93  & 4076                                                                          & 25.46  & 21374                                                                         & 130.72  & 17295                                                                         & 105.94 \\
4                     & 53290                                                                         & 778.97 & 11546                                                                         & 170.45 & 81760                                                                         & 1194.29 & 53725                                                                         & 785.30 \\ \hline \hline
\end{tabular}}
\label{Table::Cost-aComptreq}
\end{table}

The results in Table \ref{Table::Cost-aComptreq} also show that the effectiveness of action elimination depend on service capacity. The number of allowable actions varies depending on the service capacity of the system. For a given arrival volume, there is a capacity level beyond which actions cannot be eliminated based on Propositions 1 through 3. Therefore, the success of the action elimination strategy diminishes as the system capacity increase.  On the other hand, the results are independent of the customer segmentation and load patterns, since the actions are identified as suboptimal based on only the costs, $A$ and $M$. 
Overall, our results show that eliminating suboptimal policies reduces the computational requirements of the LP. For instance, for $M=2$, in each model, the number of decision variables reduces approximately to one-third, and as a result, the memory requirement decreases. 

In the experiments presented in Table \ref{ChangeCr}, we only change the rejection cost $c_r$ to 100 and 50 keeping other costs parameters the same as before. When $c_r$ decreases, the service provider can reject low-priority jobs to create capacity for high-priority jobs and avoid associated overtime and early service costs. As a result, when we compare the low-priority segmented (LS), equal segmented (ES), and high-priority segmented (HS) systems, we observe that the optimal cost is highest in HS systems and lowest in LS systems. This observation is intuitive since high-priority jobs are less flexible and costlier to serve early.

\begin{table}[h]
\caption{Computational results for $K=2$, $M=2$. }\label{ChangeCr}
\centering
\resizebox{0.9\textwidth}{!}{
\begin{tabular}{lr|rrr|rrr|rrr}
\hline \hline
                    & \multicolumn{1}{l|}{} & \multicolumn{9}{c}{$c_o=200$, $c_r=50$, $c_e^1=100$, $c_e^2=50$}                                                                                                                                                                                                    \\ \hline
                    & \multicolumn{1}{l|}{} & \multicolumn{3}{c|}{Low-priority}                                                     & \multicolumn{3}{c|}{Equal-segmented}                                                & \multicolumn{3}{c}{High-priority}                                                   \\ \hline
                    & A                     & \multicolumn{1}{c}{Cost} & \multicolumn{1}{c}{$DLPt(sec)$} & \multicolumn{1}{c|}{$RLPt(sec)$} & \multicolumn{1}{c}{Cost} & \multicolumn{1}{c}{$DLPt(sec)$} & \multicolumn{1}{c|}{$RLPt(sec)$} & \multicolumn{1}{c}{Cost} & \multicolumn{1}{c}{$DLPt(sec)$} & \multicolumn{1}{c}{$RLPt(sec)$} \\ \hline
\multirow{3}{*}{EL} & 2                     & 5.32                     & 0.66                        & 0.40                         & 8.56                     & 0.75                        & 0.36                         & 10.54                    & 0.79                        & 0.39                        \\
                    & 3                     & 17.59                    & 13.94                       & 4.15                         & 24.72                    & 14.59                       & 3.75                         & 36.47                    & 13.18                       & 4.01                        \\
                    & 4                     & 33.47                    & 198.28                      & 34.04                        & 47.88                    & 366.34                      & 48.10                        & 74.29                    & 319.37                      & 41.43                       \\ \hline
\multirow{3}{*}{FL} & 2                     & 3.74                     & 0.60                        & 0.39                         & 5.85                     & 0.61                        & 0.39                         & 6.97                     & 0.68                        & 0.38                        \\
                    & 3                     & 13.74                    & 12.91                       & 4.92                         & 21.83                    & 13.47                       & 4.81                         & 31.87                    & 13.72                       & 4.55                        \\
                    & 4                     & 28.79                    & 524.59                      & 34.53                        & 45.49                    & 194.13                      & 40.86                        & 70.08                    & 176.09                      & 39.13                       \\ \hline
\multirow{3}{*}{BL} & 2                     & 4.12                     & 0.79                        & 0.43                         & 6.79                     & 0.66                        & 0.39                         & 7.66                     & 0.83                        & 0.39                        \\
                    & 3                     & 15.11                    & 13.92                       & 4.41                         & 23.02                    & 13.92                       & 4.73                         & 30.46                    & 13.96                       & 4.59                        \\
                    & 4                     & 31.60                    & 180.91                      & 35.33                        & 47.62                    & 232.50                      & 77.72                        & 66.67                    & 270.28                      &   54.33                          \\ \hline
                    & \multicolumn{1}{l|}{} & \multicolumn{9}{c
                    }{$c_o=200$, $c_r=100$, $c_e^1=100$, $c_e^2=50$}                                                                                                                                                                                                   \\ \hline
                    & \multicolumn{1}{l|}{} & \multicolumn{3}{c|}{Low-priority}                                                     & \multicolumn{3}{c|}{Equal-segmented}                                                & \multicolumn{3}{c}{High-priority}                                                   \\ \hline
                    & A                     & \multicolumn{1}{c}{Cost} & \multicolumn{1}{c}{$DLPt(sec)$} & \multicolumn{1}{c|}{$RLPt(sec)$} & \multicolumn{1}{c}{Cost} & \multicolumn{1}{c}{$DLPt(sec)$} & \multicolumn{1}{c|}{$RLPt(sec)$} & \multicolumn{1}{c}{Cost} & \multicolumn{1}{c}{$DLPt(sec)$} & \multicolumn{1}{c}{$RLPt(sec)$} \\ \hline
\multirow{3}{*}{EL} & 2                     & 8.47                     & 0.67                        & 0.55                         & 10.83                    & 0.70                        & 0.51                         & 11.84                    & 0.68                        & 0.53                        \\
                    & 3                     & 26.68                    & 13.44                       & 5.70                         & 33.32                    & 15.17                       & 6.13                         & 41.08                    & 13.76                       & 5.74                        \\
                    & 4                     & 54.64                    & 319.59                      & \multicolumn{1}{l|}{191.39}  & 66.67                    & 206.92                      & 64.83                        & 83.23                    & 314.19                      & 240.99                      \\ \hline
\multirow{3}{*}{FL} & 2                     & 6.34                     & 0.65                        & 0.51                         & 8.91                     & 0.58                        & 0.52                         & 8.67                     & 0.66                        & 0.48                        \\
                    & 3                     & 24.56                    & 12.67                       & 5.18                         & 31.69                    & 14.63                       & 6.60                         & 36.69                    & 12.71                       & 5.40                        \\
                    & 4                     & 52.22                    & 156.39                      & 50.86                        & 64.35                    & 181.67                      & 74.48                        & 80.21                    & 173.76                      & 50.27                       \\ \hline
\multirow{3}{*}{BL} & 2                     & 6.05                     & 0.69                        & 0.56                         & 8.44                     & 0.65                        & 0.56                         & 8.41                     & 0.69                        & 0.52                        \\
                    & 3                     & 22.80                    & 13.72                       & 5.22                         & 29.53                    & 14.54                       & 7.30                         & 33.43                    & 14.48                       & 7.14                        \\
                    & 4                     & 50.12                    & 183.72                      & 66.94                        & 61.66                    & 205.09                      & 115.85                       & 74.80                    & 249.12                      & 88.94                       \\ \hline \hline
\end{tabular}}
\label{Table::Cost-12}
\end {table}

We observe a different relationship when the rejection cost increases compared to early service costs. As $c_r$ increases, it becomes less preferable to reject a low-priority job.  The results show that an equal segmented system can be costlier than a high-priority segmented system. This is especially the case when $c_r$ increases (e.g., for $c_r=100$ or $c_r=150$). However, if the arrival volume increase, we see that overtime and early service becomes more frequent, and the high-priority segmented system becomes the costliest again. In Appendix \ref{Ap:Compt}, we provide additional computational results on different cost parameters and discuss the results. These results show that when the early service costs are low compared to other cost parameters, action elimination works well and the number of actions that need to be considered is reduced by about 30\% (see Tables \ref{Table::Cost-aComptreq} and \ref{Table::6}). 

\begin{table}[h]
\caption{Computational results for $K=2$, $M=2$, $c_o=200$, $c_r=150$, $c_e^1=200$, $c_e^2=50$ }
\centering
\resizebox{0.9\textwidth}{!}{
\begin{tabular}{lr|rrr|rrr|rrr}
\hline \hline
                    & \multicolumn{1}{l|}{} & \multicolumn{3}{c|}{Low-priority}                                                     & \multicolumn{3}{c|}{Equal-segmented}                                                & \multicolumn{3}{c}{High-priority}                                                   \\ \hline
                    & A                     & \multicolumn{1}{c}{Cost} & \multicolumn{1}{c}{$DLPt(sec)$} & \multicolumn{1}{c|}{$RLPt(sec)$} & \multicolumn{1}{c}{Cost} & \multicolumn{1}{c}{$DLPt(sec)$} & \multicolumn{1}{c|}{$RLPt(sec)$} & \multicolumn{1}{c}{Cost} & \multicolumn{1}{c}{$DLPt(sec)$} & \multicolumn{1}{c}{$RLPt(sec)$} \\ \hline
\multirow{3}{*}{EL} & 2                     & 10.83                    & 0.71                        & 0.54                         & 13.05                    & 0.74                        & 0.47                         & 13.01                    & 0.70                        & 0.53                        \\
                    & 3                     & 35.54                    & 13.94                       & 6.30                         & 40.84                    & 15.55                       & 5.78                         & 45.64                    & 14.57                       & 5.68                        \\
                    & 4                     & 74.40                    & 199.17                      & 62.04                        & 82.91                    & 222.81                      & 54.75                        & 94.28                    & 213.89                      & 49.68                       \\ \hline
\multirow{3}{*}{FL} & 2                     & 8.93                     & 0.63                        & 0.47                         & 11.89                    & 0.65                        & 0.53                         & 10.31                    & 0.73                        & 0.49                        \\
                    & 3                     & 35.00                    & 12.03                       & 5.45                         & 40.98                    & 15.35                       & 6.33                         & 41.93                    & 12.85                       & 4.93                        \\
                    & 4                     & 73.07                    & 378.12                      & 69.58                        & 83.30                    & 453.90                      & 79.45                        & 90.40                    & 184.50                      & 92.74                       \\ \hline
\multirow{3}{*}{BL} & 2                     & 7.98                     & 0.72                        & 0.51                         & 10.09                    & 0.69                        & 0.56                         & 9.17                     & 0.81                        & 0.51                        \\
                    & 3                     & 30.53                    & 14.34                       & 6.28                         & 36.34                    & 14.50                       & 6.83                         & 39.27                    & 15.00                       & 6.54                        \\
                    & 4                     & 66.21                    & 192.32                      & 78.21                        & 76.49                    & 241.20                      & 140.69                       & 86.50                    & 280.27                      & 77.51                       \\ \hline \hline
\end{tabular}}
\label{Table::Cost-3}
\end{table}

Among all parameters, the length of the customer order horizon $K$ emerges as the most influential parameter on the complexity of the problem. We present the computational results for larger customer planning horizons in Table \ref{Table::Kincreases} for $c_o=200$, $c_r=150$, $c_e^1=300$ and $c_e^2=250$. The computational challenge imposed by longer customer order horizons results in longer computation time and larger memory requirements for the LP. The computation time increases rapidly as $K$ increases; however, we see that the action elimination significantly reduces running time.

\begin{table}[h]
\caption{Computational results on equal-segmented systems with $A=1$, $M=2$.}
\centering
\resizebox{0.9\textwidth}{!}{
\begin{tabular}{r|rrr|rrr|rrr}
\hline \hline
\multicolumn{1}{l|}{} & \multicolumn{3}{c|}{Back-loaded}                                                      & \multicolumn{3}{c|}{Equal-loaded}                                                     & \multicolumn{3}{c}{Front-loaded}                                                    \\ \hline
$K$                   & \multicolumn{1}{c}{Cost} & \multicolumn{1}{c}{$DLPt(sec)$} & \multicolumn{1}{c|}{$RLPt(sec)$} & \multicolumn{1}{c}{Cost} & \multicolumn{1}{c}{$DLPt(sec)$} & \multicolumn{1}{c|}{$RLPt(sec)$} & \multicolumn{1}{c}{Cost} & \multicolumn{1}{c}{$DLPt(sec)$} & \multicolumn{1}{c}{$RLPt(sec)$} \\ \hline
2                     & 0.56                     & 0.24                        & 0.02                         & 0.85                     & 0.04                        & 0.02                         & 0.51                     & 0.05                        & 0.04                        \\
3                     & 0.95                     & 9.53                        & 2.66                         & 1.36                     & 8.91                        & 2.21                         & 0.90                     & 8.38                        & 2.28                        \\
4                     & NA                       & $>7200$                     & $>7200$                      & 1.67                     & $>7200$                     & 6541.84                      & NA                       & $>7200$                     & $>7200$                     \\ \hline \hline
\end{tabular}}
\label{Table::Kincreases}
\end{table}

\subsection{Comparison of the Aggregate MDP Model to the Original MDP Model}

We perform a second set of experiments to evaluate the accuracy of of the results obtained using the proposed aggregate MDP. We compare the policies found by solving our aggregate MDP model with the policies that are optimal for the original problem. We also include in our comparison the policies found by a benchmark state aggregation method; i.e., $\varepsilon$-homogeneity (\cite{even2003approximate}, \cite{Dean2013}). In $\varepsilon$-homogeneity aggregation method, equivalence is measured in terms of one-step transition probabilities and costs and $\varepsilon$-equivalent states are aggregated. Level of approximation provided by the aggregate MDP depends on $\varepsilon \in [ 0, 1 ]$ value chosen. 

To compare aggregate MDP models to the original model and with each other, we define two performance indicators: absolute percentage gap (AG) and matched action percentage (MAP). MAP is the percentage of states with identical actions under the two policies compared out of all states in the state space (i.e., all $|S|$ states).
AG demonstrates the percent deviation from the optimal long-run average cost. MAP shows the average accuracy of the aggregate MDP model in terms of its ability to find the optimal decision across all states of the original MDP model. 

In Table \ref{Table::Aggregationsmall-woAE}, we use $t$ to denote the total computation time and $Agg-t$ to denote the time needed to perform aggregation, both measured in seconds. The optimal solution is indicated as OPT, and approximate solutions are named after the state aggregation method, e.g., Tot-job, $\varepsilon$-agg. We use an equal-segmented system with $K=2$, $A=5$ and $M=2$ as reference. We present the computational results for $c_o=200$, $c_r=150$, $c_e^1=100$ and $c_e^2=50$ and three load types. The results in Table \ref{Table::Aggregationsmall-woAE} are obtained without action elimination. Approximate optimal policies are evaluated according to equation \eqref{Aggregation::Evaluation}, and feasible policies are obtained for the original MDP.

These results show that our aggregation method successfully approximates the MDP model and is a powerful alternative to $\varepsilon$-homogeneity, even though it generates relatively larger aggregate MDP models. The total-job aggregation method results in more accurate aggregate MDP models with AG ranging from 2\% to 15\%. In terms of computation time, our aggregation method outperforms $\varepsilon$-homogeneity approximations. We also implemented the total-job aggregation along with action elimination. However, due to the elimination of common actions for states, the level of aggregation is reduced. As a result, the efficiency of the aggregation decreased.

\begin{table}[h]
   \caption{Computational results for an equal segmented system with $K=2$, $A=5$, $M=2$ without action elimination.}
    \begin{adjustwidth}{-.5in}{-.5in}  
\centering
\resizebox{1.1\textwidth}{!}{
\begin{tabular}{r|rrrrrrrr|rrrrrrrr|rrrrrrrr}
\hline \hline
\multicolumn{1}{r|}{} & \multicolumn{8}{c|}{Back-loaded} & \multicolumn{8}{c|}{Equal-loaded}  & \multicolumn{8}{c}{Front-loaded}         \\ \hline
                      & \makecell{Number of \\states} & \makecell{Number of \\actions} & Cost   & $Mem$  & $t$     & $Agg-t$ & MAP   & AG    & \makecell{Number of \\states} & \makecell{Number of \\actions} & Cost   & $Mem$  & $t$     & $Agg-t$ & MAP   & AG    & \makecell{Number of \\states} & \makecell{Number of \\actions} & Cost   & $Mem$  & $t$    & $Agg-t$ & MAP   & AG   \\ \hline
OPT                   & 14256                                                                        & 180435                                                                        & 129.02 & 815.14 & 1745.10 & NA      & NA    & NA    & 14256                                                                        & 180435                                                                        & 135.11 & 815.14 & 1345.54 & NA      & NA    & NA    & 14256                                                                        & 180435                                                                        & 136.44 & 815.14 & 827.01 & NA      & NA    & NA   \\
0.2-agg               & 623                                                                          & 4588                                                                          & 155.86 & 81.06  & 844.80  & 291.87  & 45.86 & 20.81 & 616                                                                          & 4622                                                                          & 147.91 & 81.13  & 863.17  & 279.19  & 45.92 & 9.48  & 619                                                                          & 4602                                                                          & 141.07 & 82.44  & 851.46 & 275.87  & 54.67 & 3.40 \\
0.4-agg               & 615                                                                          & 4828                                                                          & 156.90 & 86.36  & 859.28  & 278.65  & 45.82 & 21.62 & 613                                                                          & 4576                                                                          & 147.80 & 80.41  & 841.20  & 274.49  & 45.77 & 9.39  & 622                                                                          & 4976                                                                          & 140.02 & 88.05  & 883.83 & 283.86  & 54.79 & 2.63 \\
0.6-agg               & 617                                                                          & 4683                                                                          & 158.01 & 82.91  & 834.53  & 282.82  & 45.64 & 22.47 & 616                                                                          & 4676                                                                          & 148.26 & 82.03  & 825.55  & 274.43  & 45.81 & 9.74  & 613                                                                          & 4615                                                                          & 140.81 & 81.24  & 868.64 & 291.06  & 54.93 & 3.21 \\
0.8-agg               & 621                                                                          & 4919                                                                          & 157.85 & 86.27  & 892.59  & 276.30  & 45.79 & 22.35 & 614                                                                          & 4590                                                                          & 150.00 & 81.93  & 855.47  & 284.79  & 45.78 & 11.02 & 620                                                                          & 4832                                                                          & 142.09 & 86.47  & 870.23 & 305.21  & 54.84 & 4.15 \\
1-agg                 & 618                                                                          & 4880                                                                          & 159.83 & 86.46  & 854.67  & 272.71  & 45.89 & 23.88 & 620                                                                          & 4686                                                                          & 151.27 & 82.65  & 832.96  & 275.68  & 45.76 & 11.96 & 616                                                                          & 4729                                                                          & 142.15 & 83.39  & 857.91 & 275.17  & 54.73 & 4.19 \\
Tot-job               & 1673                                                                         & 8983                                                                          & 148.47 & 355.64 & 448.77  & 63.80   & 45.72 & 15.08 & 1875                                                                         & 10396                                                                         & 143.00 & 477.57 & 486.31  & 62.09   & 45.75 & 5.84  & 1567                                                                         & 8166                                                                          & 139.21 & 312.87 & 473.78 & 68.74   & 54.67 & 2.03 \\ \hline \hline
\end{tabular}}
\label{Table::Aggregationsmall-woAE}
    \end{adjustwidth}

\end{table}

Note that an approximate optimal policy can be easily computed even for very large instances. However, evaluating the resulting approximate optimal policy back in the original MDP is challenging and sometimes intractable, as it also requires the use of the full state space. In such cases, 
the approximate optimal policy provides insight about the behavior of the solution.

\subsection{Comparison of the Aggregate MDP Model to Benchmarks}

In the third set of computational experiments, we perform a simulation study on larger instances where the LP formulation cannot be solved due to memory limits. We compare our aggregate MDP model (i.e., Tot-job) with two simple benchmark heuristics. We also report results of $\varepsilon$-aggregate models. The first benchmark heuristic is a myopic policy (MP) that picks the minimum-cost action at each state. The second benchmark policy, which we call always-serve policy (AP), forces the system to work in full capacity and aims to serve as many jobs as possible. 

In the simulation study, we perform 1.1 million iterations for each policy for each parameter setting (200K for warm-up and 900K for evaluations). This study shows the trade-off between the solution quality and the computation time and memory requirements.

We run the simulations using the following cost parameters: $c_o=200$, $c_r=150$, $c_e^1=100$, $c_e^2=50$. Table \ref{Table::Simulation1} shows the results for an equal segmented system with $K=2$ and $A=6$, $M=2$. We compare approximate and heuristic results for this instance under different customer order patterns. In Table \ref{Table::Simulation1}, $t$ refers to the time until a policy is found. We do not report the duration of simulation since it is constant once a policy is determined. 

\begin{table}[H]
\caption{Simulation results for an equal segmented system with $K=2$, $A=6$, $M=2$.}
    \begin{adjustwidth}{-.5in}{-.5in}  

\centering
\resizebox{1.1\textwidth}{!}{
\begin{tabular}{r|rrrrrr|rrrrrr|rrrrrr}
\hline \hline
\multicolumn{1}{l|}{} & \multicolumn{6}{c|}{Back-loaded}                                                                                                                                                                 & \multicolumn{6}{c|}{Equal-loaded}                                                                                                                                                               & \multicolumn{6}{c}{Front-loaded}                                                                                                                                                                 \\ \hline
                      & \makecell{Number of \\states} & \makecell{Number of \\actions} & Cost   & $Mem$  & $t$    & $Agg-t$ & \makecell{Number of \\states} & \makecell{Number of \\actions} & Cost    & $Mem$  & $t$    & $Agg-t$ & \makecell{Number of \\states} & \makecell{Number of \\actions} & Cost    & $Mem$  & $t$    & $Agg-t$ \\ \hline
0.5-agg               & 943                                                                         & 1038                                                                          & 496.88 & 25.52  & 429.56 & 0.66    & 943                                                                         & 1038                                                                          & 547.68  & 25.52  & 424.63 & 0.86    & 943                                                                         & 1038                                                                          & 621.55  & 25.52  & 414.32 & 0.65    \\
1-agg                 & 920                                                                         & 938                                                                           & 538.82 & 22.73  & 401.43 & 0.70    & 920                                                                         & 938                                                                           & 597.25  & 22.73  & 402.55 & 0.74    & 920                                                                         & 938                                                                           & 680.98  & 22.73  & 405.53 & 0.65    \\
Tot-job               & 4868                                                                        & 5317                                                                          & 86.53  & 558.05 & 620.03 & 0.32    & 4868                                                                        & 5317                                                                          & 265.60  & 558.05 & 588.55 & 0.34    & 4868                                                                        & 5317                                                                          & 550.00  & 558.05 & 608.26 & 0.32    \\
MP                    & 31213                                                                       & 510510                                                                        & 413.58 & 1.93   & 51.95  & NA      & 31213                                                                       & 510510                                                                        & 603.96  & 1.93   & 52.58  & NA      & 31213                                                                       & 510510                                                                        & 2843.00 & 1.93   & 52.65  & NA      \\
AP                    & 31213                                                                       & 510510                                                                        & 390.82 & 1.93   & 51.95  & NA      & 31213                                                                       & 510510                                                                        & 1112.83 & 1.93   & 52.59  & NA      & 31213                                                                       & 510510                                                                        & 2918.66 & 1.93   & 52.66  & NA      \\ \hline \hline
\end{tabular}}
\end{adjustwidth}
\label{Table::Simulation1}
\end{table}

Aggregate MDP models require higher computational efforts compared to heuristic policies. In terms of computation time, the heuristic policies beats the performance of aggregate MDP models in every setting. On the other hand, the aggregate MDP models provide much better results in terms of long-run average cost in almost every setting.

Based on our computational experiments, we observe that the service provider's optimal strategy strongly depends on the service capacity and operational costs. For small systems with short customer order horizon and small customer arrivals, once a set of suboptimal policies are excluded from the policy space through action elimination, we can find the optimal policy within a reasonable time, e.g., in seconds. Since the optimal policy is a deterministic one, the service provider can directly use it until one of the problem parameters changes. An increase in server capacity or an increase in maximum number of arrivals, or a change in customer order pattern do not substantially affect the performance of the solution method. For large instances, when the quality of the solution is more critical than the computation time, using Total-job aggregate MDP model without action elimination is the best solution method. However, simple heuristics, i.e., MP and AP, can also provide good solutions depending on the server capacity and costs. Hence, simple heuristics can also be used when the computation time matters.

\section{Conclusion}
\label{Sec:Conclusion}
In this paper, we study admission and capacity allocation problem in a queuing system where customers have different priorities and submit their desired service delivery time before joining the system. 
In our system, customers' perception of the service depends on the service provider's ability to serve the service requests as close as possible to their preferred delivery times. Therefore, any deviation is penalized. Service provider's problem is determine which arriving jobs should be admitted and when to serve each job considering the preferred service time and the priority of customers with minimum cost. 

We formulate the service provider's problem as a discrete-time finite-horizon MDP, and develop solution methodologies based on linear programming. The size of the MDP model increases exponentially based on the length of the customer order horizon. To increase efficiency, we partially characterize a set of suboptimal policies. Accordingly, we decrease the problem size and the required computation time to find optimal solutions significantly. In addition, we develop a state-space aggregation model based on our new \emph{equivalence} definition to obtain heuristic solutions. We perform extensive computational experiments to compare our heuristic with other other aggregations and benchmark heuristics. On average, our method outperforms others in terms of solution quality especially significantly on larger instances.

There are several future research directions. We assume that each job can be completed in one period, and arrival rate is constant. Considering more general job completions and customer arrivals based on popular times are potential directions. In addition, it might be of interest to consider a game-theoretical setting where customers determine their submitted preferred times based on the information on the service provider's historical performance.


\begin{singlespacing}
\small
\FloatBarrier
\bibliographystyle{chicago}
 
\bibliography{thesislib}

\begin{thebibliography}{}

\bibitem[\protect\citeauthoryear{Aho}{Aho}{2019}]{Curbside}
Aho, C. (2019).
\newblock Curbside pickup is taking over the market.
\newblock
  \url{https://usa.inquirer.net/35574/curbside-pickup-is-taking-over-the-market}.
\newblock Accessed: July 3, 2020.

\bibitem[\protect\citeauthoryear{Ayvaz and Huh}{Ayvaz and
  Huh}{2010}]{Ayvaz2010}
Ayvaz, N. and W.~T. Huh (2010).
\newblock {Allocation of hospital capacity to multiple types of patients}.
\newblock {\em Journal of Revenue and Pricing Management\/}~{\em 9\/}(5),
  386--398.

\bibitem[\protect\citeauthoryear{Bhulai and Koole}{Bhulai and
  Koole}{2003}]{Bhulai}
Bhulai, S. and G.~Koole (2003).
\newblock {On the Structure of Value Functions for Threshold Policies in
  Queueing Models}.
\newblock {\em {Journal of Applied Probability}\/}~{\em 40\/}(3), 613--622.

\bibitem[\protect\citeauthoryear{Carr and Duenyas}{Carr and
  Duenyas}{2000}]{Carr2000}
Carr, S. and I.~Duenyas (2000).
\newblock {Optimal admission control and sequencing in a
  make-to-stock/make-to-order production system}.
\newblock {\em Operations Research\/}~{\em 48\/}(5), 709--720.

\bibitem[\protect\citeauthoryear{{\c{C}}avdar and I\c{s}{\i}k}{{\c{C}}avdar and
  I\c{s}{\i}k}{2022}]{Cavdar2020}
{\c{C}}avdar, B. and T.~I\c{s}{\i}k (2022).
\newblock Capacity allocation in queuing systems with preferred service
  completion times.
\newblock {\em Naval Research Logistics (NRL)\/}, 1--31.
  https://doi.org/10.1002/nav.22046.

\bibitem[\protect\citeauthoryear{Chopra}{Chopra}{2016}]{Chopra2016}
Chopra, S. (2016).
\newblock {How omni-channel can be the future of retailing}.
\newblock {\em Decision\/}~{\em 43\/}(2), 135--144.

\bibitem[\protect\citeauthoryear{Company}{Company}{2015}]{Nielsen2015}
Company, T.~N. (2015).
\newblock {The Future of Grocery E- Commerce}.
\newblock ~(April), 1--35.

\bibitem[\protect\citeauthoryear{Danziger}{Danziger}{2019}]{Forbes}
Danziger, P.~N. (2019).
\newblock Walmart leads the soon-to-be \$35 billion curbside pickup market.
\newblock
  \url{https://www.forbes.com/sites/pamdanziger/2019/04/07/walmart-is-in-the-lead-in-the-soon-to-be-35-billion-\Urlbreak
  curbside-pickup-market/#6a19d996199e}.
\newblock Accessed: August 10, 2020.

\bibitem[\protect\citeauthoryear{Dean, Givan, and Leach}{Dean
  et~al.}{2013}]{Dean2013}
Dean, T.~L., R.~Givan, and S.~Leach (2013).
\newblock {Model Reduction Techniques for Computing Approximately Optimal
  Solutions for Markov Decision Processes}.
\newblock {\em arXiv preprint arXiv:1302.1533\/}, 124--131.

\bibitem[\protect\citeauthoryear{Duenyas}{Duenyas}{1995}]{Duenyas1995}
Duenyas, I. (1995).
\newblock {Single Facility Due Date Setting with Multiple Customer Classes}.
\newblock {\em Management Science\/}~{\em 41\/}(4), 608--619.

\bibitem[\protect\citeauthoryear{Even-Dar and Mansour}{Even-Dar and
  Mansour}{2003}]{even2003approximate}
Even-Dar, E. and Y.~Mansour (2003).
\newblock Approximate equivalence of markov decision processes.
\newblock In {\em Learning Theory and Kernel Machines}, pp.\  581--594.
  Springer.

\bibitem[\protect\citeauthoryear{Givan, Dean, and Greig}{Givan
  et~al.}{2003}]{Givan2003}
Givan, R., T.~Dean, and M.~Greig (2003).
\newblock {Equivalence notions and model minimization in Markov decision
  processes}.
\newblock {\em Artificial Intelligence\/}~{\em 147\/}(1-2), 163--223.

\bibitem[\protect\citeauthoryear{Gocgun and Puterman}{Gocgun and
  Puterman}{2014}]{Gocgun2014}
Gocgun, Y. and M.~L. Puterman (2014).
\newblock {Dynamic scheduling with due dates and time windows: An application
  to chemotherapy patient appointment booking}.
\newblock {\em Health Care Management Science\/}~{\em 17\/}(1), 60--76.

\bibitem[\protect\citeauthoryear{Green, Savin, and Wang}{Green
  et~al.}{2006}]{Green2006}
Green, L.~V., S.~Savin, and B.~Wang (2006).
\newblock {Managing patient service in a diagnostic medical facility}.
\newblock {\em Operations Research\/}~{\em 54\/}(1), 11--25.

\bibitem[\protect\citeauthoryear{Kolesar}{Kolesar}{1970}]{Kolesar1970}
Kolesar, P. (1970).
\newblock {A Markovian Model for Hospital Admission Scheduling}.
\newblock {\em Management Science\/}~{\em 16\/}(6), B--384--B--396.

\bibitem[\protect\citeauthoryear{Lasserre}{Lasserre}{1994}]{Lasserre1994}
Lasserre, J.~B. (1994).
\newblock {Detecting Optimal and Non-Optimal Actions in Average-Cost Markov
  Decision Processes}.
\newblock {\em Journal of Applied Probability\/}~{\em 31\/}(4), 979--990.

\bibitem[\protect\citeauthoryear{Lewis and Slotnick}{Lewis and
  Slotnick}{2002}]{Lewis2002}
Lewis, H.~F. and S.~A. Slotnick (2002).
\newblock {Multi-period job selection: Planning work loads to maximize profit}.
\newblock {\em Computers \& Operations Research\/}~{\em 29\/}(8), 1081--1098.

\bibitem[\protect\citeauthoryear{Papadimitriou and Tsitsiklis}{Papadimitriou
  and Tsitsiklis}{1987}]{Papadimitriou1987}
Papadimitriou, C.~H. and J.~N. Tsitsiklis (1987).
\newblock {The Complexity of Markov Decision Processes}.
\newblock {\em Mathematics of Operations Research\/}~{\em 12\/}(3), 441--450.

\bibitem[\protect\citeauthoryear{Patrick, Puterman, and Queyranne}{Patrick
  et~al.}{2008}]{Patrick2008}
Patrick, J., M.~L. Puterman, and M.~Queyranne (2008).
\newblock {Dynamic multipriority patient scheduling for a diagnostic resource}.
\newblock {\em Operations Research\/}~{\em 56\/}(6), 1507--1525.

\bibitem[\protect\citeauthoryear{Puterman}{Puterman}{1994}]{Puterman}
Puterman, M.~L. (1994).
\newblock {\em Markov Decision Processes: Discrete Stochastic Dynamic
  Programming}.
\newblock New York: John Wiley \& Sons.

\bibitem[\protect\citeauthoryear{Sheehan and Skelly}{Sheehan and
  Skelly}{2021}]{Sheehan2021}
Sheehan, E. and L.~Skelly (2021).
\newblock {The Click and Collect Consumer}.
\newblock {\em Deloitte\/}.

\bibitem[\protect\citeauthoryear{Slotnick}{Slotnick}{2011}]{Slotnick2011}
Slotnick, S.~A. (2011).
\newblock {Order acceptance and scheduling: A taxonomy and review}.
\newblock {\em European Journal of Operational Research\/}~{\em 212\/}(1),
  1--11.

\bibitem[\protect\citeauthoryear{Slotnick and Morton}{Slotnick and
  Morton}{1996}]{Slotnick1996}
Slotnick, S.~A. and T.~E. Morton (1996).
\newblock Selecting jobs for a heavily loaded shop with lateness penalties.
\newblock {\em Computers \& Operations Research\/}~{\em 23\/}(2), 131--140.

\bibitem[\protect\citeauthoryear{Somers}{Somers}{2020}]{Somers2020}
Somers, A.~R. (2020).
\newblock {The consumer transformed - Changing behaviours are accelerating
  trends along a reinvented customer purchase journey}.
\newblock {\em PwC\/}~(May), 1--25.

\bibitem[\protect\citeauthoryear{Standish and Bossi}{Standish and
  Bossi}{2020}]{Standish2020}
Standish, J. and M.~Bossi (2020).
\newblock {How is COVID-19 changing the retail consumer? - Data-driven insights
  into consumer behavior}.
\newblock {\em Accenture\/}~(August), 1--32.

\bibitem[\protect\citeauthoryear{Sutton and Barto}{Sutton and
  Barto}{2018}]{sutton2018reinforcement}
Sutton, R.~S. and A.~G. Barto (2018).
\newblock {\em Reinforcement learning: An introduction}.
\newblock MIT press.

\end{thebibliography}
\newpage
\end{singlespacing}
\appendix

\section{Flow of the service process}
\label{App:Flow}

\begin{figure}[H]
\includegraphics[scale=0.6]{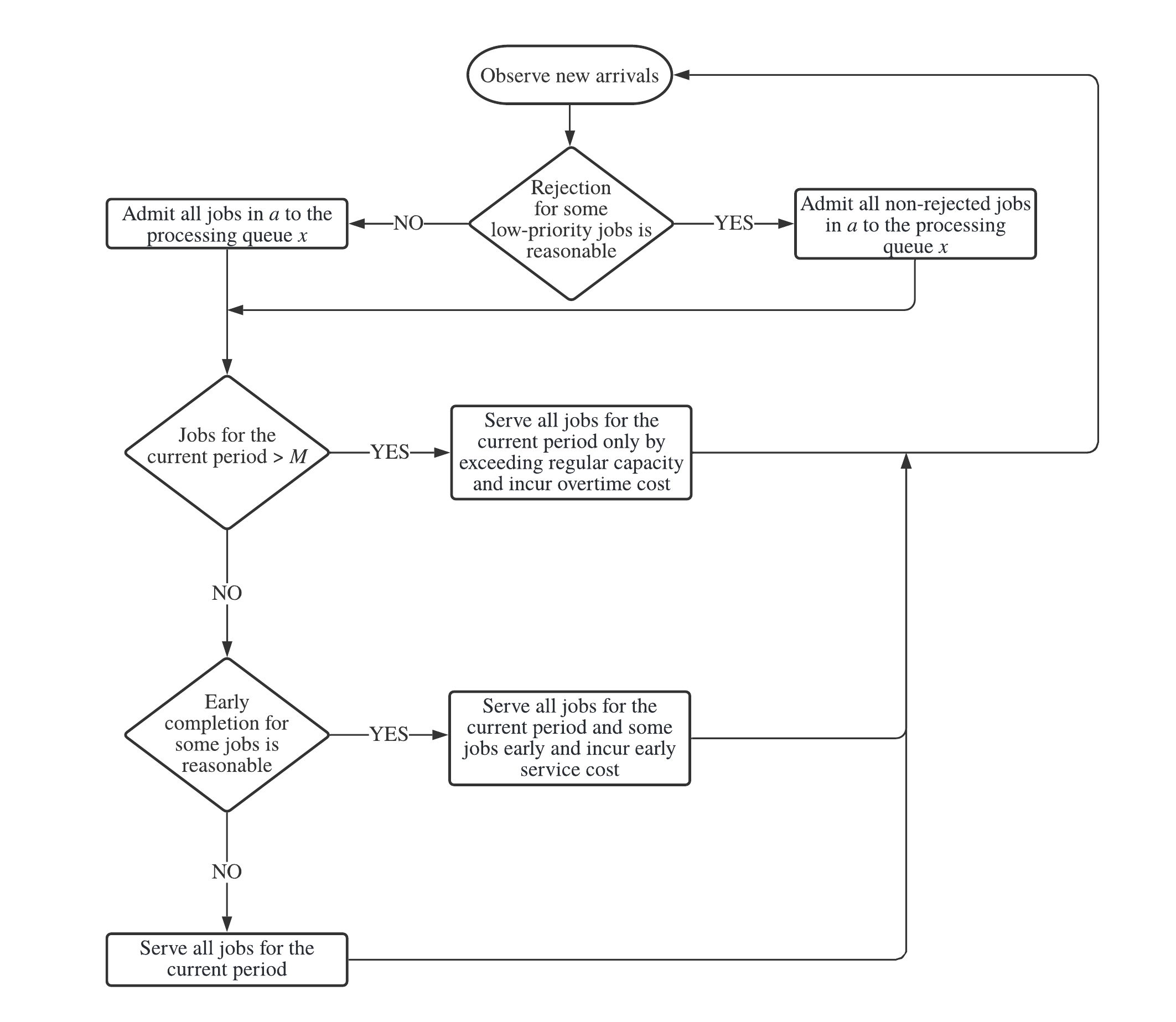}
\centering
\caption{Flowchart of the service process.}
\label{Fig:Flow}

\end{figure}

\section{Proposition Proofs}
\label{Ap:Pros}

\subsection{Proof of Proposition \ref{Proposition::1}}

\begin{proof}
Suppose that a high-priority job that is due in period $t + j$ is served early in overtime. Then, a cost of $j c_{1}^e + c_o$ is incurred for this job. However, if this job was kept in the queue until its requested delivery time, it would be served in overtime with some probability $p, 0 \leq p \leq 1$, at the cost expected cost of $p c_o$. Since $j c_{1}^e + c_o > p c_o$, a job should never be early served in overtime. An analogous argument is valid for low-priority jobs. 
\end{proof}

\subsection{Proof of Proposition \ref{Proposition::2}}

\begin{proof}
We start with the rejection decision of low-priority jobs for the current period, i.e., $min\{ [(x_{10}+a_{10}+a_{20})-M]^+,a_{20} \}  \le r_{20}\le a_{20}$, which is common for each cost structure. After observing the acceptance queue, if no rejection is made, $[(x_{10}+a_{10}+a_{20})-M]^+$ is the number of jobs that need to be fulfilled in overtime at a cost of $c_o[(x_{10}+a_{10}+a_{20})-M]^+$. However, if one of the excess low-priority jobs is rejected, the corresponding cost is $c_r+ c_o([(x_{10}+a_{10}+a_{20})-M]^+-1)$ and it is always less than the total overtime cost due to Assumption \ref{Assumption::1}. Therefore, whenever it is certain that the system operates in overtime, the excess low-priority jobs in the acceptance queue should be rejected; the minimum number of necessary rejections determines the lower bound for $r_{20}$. The upper bound on $r_{20}$ follows by definition.

\begin{itemize}
    \item When $c_{2}^e < c_{1}^e < c_{r} < c_{o}$ : 
    Eq. $(8)$ suggests that keeping all low-priority requests to serve on time is suboptimal since, their rejection or early service would be less costly. $a_{1 1} + a_{2 1}$ is the number of arriving request for $t+1$. If all is accepted, in the next period all would have to be served. Assume that $a_{1 1} + a_{2 1} > M$, it is certain that there will be a need for overtime in the next period. The overtime cost can be avoided by rejecting or early serving some excess low-priority jobs. Due to Proposition \ref{Proposition::1}, we can only early serve if there are some idle server capacity. Therefore, rejection and early service decisions are interdependent, at least $min\{[[(a_{1 1}+a_{2 1}\ )-M]^+-y_{2 1} ]^+,a_{2 1}  \}$ and at most $[a_{2 1}-[M-(x_{1 0}+a_{1 0}+a_{2 0} - r_{2 0}\ )]^+ ]^+$ among low-priority job requests should be rejected. The upper bound on the service decision (i.e., Eq. $(10)$) is also determined accordingly. Similarly, in the presence of idle server capacity, keeping all high-priority requests to serve on time is suboptimal since their early service would be less costly. Therefore, we should serve some orders early to avoid the overtime cost in the next period, and the upper bound follows by the Lemma.
    
    \item When $ c_{r} < c_{2}^e < c_{1}^e < c_{o} $ : 
    Considering low-priority jobs for the next period, whenever the total number of jobs for the next period exceeds the server capacity, keeping excess jobs in the processing queue to serve on time or early are both suboptimal. When $a_{1 1} + a_{2 1} > M$, the excess low-priority job requests should be rejected at a cost of $c_r r_{21}$ in the current period instead of $c_o r_{2 1}$ in the next period, $c_r r_{21} < c_o r_{2 1}$ . Then, at least $min\{[(a_{11}+a_{21})-M]^+,a_{21}\}$ among low-priority job requests for $t+1$ should be rejected. The upper bound on the early service decision of high-priority jobs is easy to show. The idle capacity is $[M-(x_{1 0} + a_{1 0} + a_{2 0} - r_{2 0} + y_{2 1})]^+$. We have $x_{1 1}=0$ and the system can only serve newly arriving requests, up to the amount of the idle capacity.
    
    \item When $ c_{r} < c_{o} < c_{2}^e < c_{1}^e $ : 
     Let $a_{1 1} + (a_{2 1}-r_{2 1})$ be the total number of arriving jobs for $t+1$, and it is certain that $a_{1 1} + (a_{2 1}-r_{2 1}) \leq M$. For any realization of acceptance queue, in the presence of idle capacity, if one of these high-priority job requests is served in advance, the system incurs $c_e^1$ amount of early service cost. After observing the acceptance queue, assume that $a_{11} + (a_{21}-r_{21}) - 1 + a_{10}^{t+1} + (a_{20}^{t+1}-r_{20}^{t+1}) > M $, and the system incurs $c_o[(a_{11} + (a_{21}-r_{21}) - 1 + a_{10}^{t+1} + (a_{20}^{t+1}-r_{20}^{t+1}) )-M]^+ $. If the system does not serve any request early, at time $t+1$ we incur $c_o[(a_{11} + (a_{21}-r_{21}) + a_{10}^{t+1} + (a_{20}^{t+1}-r_{20}^{t+1}) )-M]^+ $, which is clearly less than $c_e^1 + c_o[(a_{11} + (a_{21}-r_{21}) - 1 + a_{10}^{t+1} + (a_{20}^{t+1}-r_{20}^{t+1}) )-M]^+ $. The observation still holds when $a_{11} + (a_{21}-r_{21}) - 1 + a_{10}^{t+1} + (a_{20}^{t+1}-r_{20}^{t+1}) < M $ after the arrivals. Therefore, even though there is some idle capacity, early service is suboptimal. Early service decision for low-priority jobs follows by the same argument.
\end{itemize}
\end{proof}

\subsection{Proof of Proposition \ref{Proposition::3}} 

\begin{proof}
The proposition is based on a simple argument. For the first part, if a state becomes $j$-boundary for all $j$ after the arrivals, then it is certain that for the next $K$ periods, the jobs that are due in that  period will exceed the capacity, and therefore there will be overtime for sure. If we reject the jobs that would cause overtime, we would incur the rejection cost as opposed to a larger overtime cost that will be incurred with probability 1. Therefore, when there is certainty that a job will be served in overtime unless it is rejected upon arrival, then this job should be rejected. The second part follows from Proposition \ref{Proposition::1}.
\end{proof}

\subsection{Proof of Proposition \ref{Proposition::4}}

\begin{proof}
The same argument in Proposition \ref{Proposition::2} to prove the rejection decision of low-priority jobs for the current period always holds. 

\begin{itemize}
    \item When $c_{2}^e < c_{1}^e < c_{r} < c_{o}$ : $a_{1 j} + a_{2 j}$ is the number of arriving jobs due $t+1$, and the total number of jobs for $t+j$ is $x_{1 j} + x_{2 j} + a_{1 j} + a_{2 j}$. If all job requests for $t+j$ are admitted and held in the queue, in $j$ periods they will all become high-priority requests that have to be served in current period. If $x_{1 j} + x_{2 j} + a_{1 j} + a_{2 j} > M$, then it is known that there will be a need for overtime for sure at $t+j$. Therefore, to avoid this overtime cost, the excess low-priority job requests should be rejected and the system incurs $c_r r_{2 j}$ in this period instead of $c_o r_{2 j}$ in $j$ periods later. At least $min\{[(x_{1 j} + x_{2 j} + a_{1 j}+a_{2 j})-M]^+,a_{2 j}\}$ among low-priority job requests for $t+j$ should be rejected. The upper bound on $r_{2 j}$ follows by definition. The upper bound on the early service decision of high-priority jobs is easy to show. Due to Proposition \ref{Proposition::1} jobs in the system can be served at most the amount of remaining idle capacity. Similarly, the upper bound on $y_{2 j}$ is settled by definition. 
   
    \item When $ c_{r} < c_{2}^e < c_{1}^e < c_{o} $ : $a_{1j}+a_{2j}$ is the number of arriving request due $t+1$, and the total number of job request for $t+j$ is $x_{1j} + x_{2j} + a_{1j}+a_{2j} $. If all job requests for $t+j$ are admitted and held in the queue, in $j$ periods they will all become high-priority customers have to be served in current period. If $x_{1j} + x_{2j} + a_{1j}+a_{2j} > M$, then it is known that there will be a need for overtime for sure at $t+j$. Therefore, to avoid this overtime cost, the excess low- priority job requests should be rejected and the system incurs $c_r r_{2j}$ in the current period instead of $c_or_{2j}$ in $t+j$. At least $min\{[(x_{1j} + x_{2j} + a_{1j}+a_{2j})-M]^+,a_{2j}\}$ among low-priority job requests for $t+j$ should be rejected. The upper bound on $r_{2j}$ follows by definition. The upper bound on the early service decision of high-priority jobs is easy to show. Due to Proposition \ref{Proposition::1} jobs in the system can be served at most the amount of remaining idle capacity. Similarly, the upper bound on $y_{2 j}$ is settled by definition.  
    
    \item When $ c_{r} < c_{o} < c_{2}^e < c_{1}^e $ : $a_{1j} + (a_{2 j}-r_{2 j})$ is the total number of arriving requests for $t+j$, and it is certain that $ x_{1 j} + x_{2 j} + a_{1 j} + (a_{2 j}-r_{2 j}) \le M$. Independent of arriving jobs, in the presence of idle capacity, if one of these high-priority job requests is served in advance, the system incurs $j c_e^1$ amount of early service cost. For $j=1$ refer to Proposition \ref{Proposition::2}. For $j>1$ since $c_{o} < c_{1}^e$, the same argument still holds. Early service decision for low-priority jobs follows by the same argument.
    
\end{itemize}

\end{proof}

\subsection{Proof of Proposition \ref{Pro:AggMarkov}}

\begin{proof}
To prove the proposition, it is enough to show that the sum of the entries in each row is equal to 1.
\begin{equation}
    \begin{split}
        \sum_{j=1}^{N} P(S_j | S_i, d) & = \sum_{j=1}^{N} \sum_{s \in S_i} \frac{1}{|S_i|} \sum_{s'\in S_j}P(s'| s, d) , \quad \forall S_i \in \mathbf{S} \quad \forall d \in \mathcal{D}_{S_i} \\ 
        & = \sum_{s \in S_i}   \frac{1}{|S_i|} \sum_{j=1}^{N} \sum_{s'\in S_j}P(s'| s, d) , \quad \forall S_i \in \mathbf{S} \quad \forall d \in \mathcal{D}_{S_i} \\ 
        & =  \frac{1}{|S_i|} \sum_{s \in S_i} \sum_{j=1}^{N}\sum_{s'\in S_j}P(s' | s, d),\quad \forall S_i \in \mathbf{S} \quad \forall d \in \mathcal{D}_{S_i} \\
        & = \frac{1}{|S_i|} \sum_{s \in S_i} \sum_{s'\in S}P(s'|s,d), \quad \forall S_i \in \mathbf{S} \quad \forall d \in \mathcal{D}_{S_i} \\
             & = \frac{1}{|S_i|} \sum_{s \in S_i} 1, \in \mathbf{S}\\
        &=1,\quad \forall S_i \in \mathbf{S} \quad \forall d \in \mathcal{D}_{S_i} .
    \end{split}
\end{equation}
\end{proof}

\section{Additional Computational Results}
\label{Ap:Compt}

\begin{table}[H]
\setcellgapes{2.5pt}\makegapedcells
\caption{Computational requirements for $K=2$, $M=2$.}
\centering
\resizebox{0.9\textwidth}{!}{
\begin{tabular}{r|rr|rr|rr|rr}
\hline \hline
\multicolumn{1}{l|}{} & \multicolumn{2}{c|}{}                                                                  & \multicolumn{2}{c|}{\makecell{$c_o=200$, $c_r=50$,\\$c_e^1=100$, $c_e^2=50$}}                       & \multicolumn{2}{c|}{\makecell{$c_o=200$, $c_r=100$,\\$c_e^1=100$, $c_e^2=50$}}                     & \multicolumn{2}{c}{\makecell{$c_o=200$, $c_r=150$,\\$c_e^1=200$, $c_e^2=50$}}                      \\ \hline
\multicolumn{1}{l|}{} & \multicolumn{2}{c|}{DLP}                                                               & \multicolumn{2}{c|}{RLP}                                                               & \multicolumn{2}{c|}{RLP}                                                               & \multicolumn{2}{c}{RLP}                                                                \\ \hline
$A$                   & \makecell{Number of \\actions} & $Mem$  & \makecell{Number of \\actions} & $Mem$  & \makecell{Number of \\actions} & $Mem$  & \makecell{Number of \\actions} & $Mem$  \\ \hline
2                     & 1896                                                                          & 4.08   & 864                                                                           & 1.97   & 1040                                                                          & 2.33   & 1047                                                                          & 2.35   \\
3                     & 12210                                                                         & 74.93  & 3262                                                                          & 20.52  & 4076                                                                          & 25.46  & 4050                                                                          & 25.31  \\
4                     & 53290                                                                         & 778.97 & 9238                                                                          & 136.82 & 11546                                                                         & 170.45 & 11383                                                                         & 168.07 \\ \hline \hline
\end{tabular}}
\label{Table::Compt-req-123}
\end{table}

\begin{table}[H]
\caption{Computational results for $K=2$, $M=2$. }
\centering
\resizebox{0.9\textwidth}{!}{
\begin{tabular}{lr|rrr|rrr|rrr}
\hline \hline
                    & \multicolumn{1}{l|}{} & \multicolumn{9}{c}{$c_o=200$, $c_r=150$, $c_e^1=300$, $c_e^2=250$}                                                                                                                                                                                                   \\ \hline
                    & \multicolumn{1}{l|}{} & \multicolumn{3}{c|}{Equal-segmented}                                                & \multicolumn{3}{c|}{High-priority}                                                    & \multicolumn{3}{c}{Low-priority}                                                     \\ \hline
                    & A                     & \multicolumn{1}{c}{Cost} & \multicolumn{1}{c}{$DLPt(sec)$} & \multicolumn{1}{c|}{$RLPt(sec)$} & \multicolumn{1}{c}{Cost} & \multicolumn{1}{c}{$DLPt(sec)$} & \multicolumn{1}{c|}{$RLPt(sec)$} & \multicolumn{1}{c}{Cost} & \multicolumn{1}{c}{$DLPt(sec)$} & \multicolumn{1}{c}{$RLPt(sec)$} \\ \hline
\multirow{3}{*}{EL} & 2                     & 14.75                    & 0.59                        & 0.25                         & 14.02                    & 0.52                        & 0.30                         & 12.11                    & 0.49                        & 0.29                        \\
                    & 3                     & 45.80                    & 14.68                       & 3.27                         & 48.29                    & 14.30                       & 3.33                         & 40.85                    & 13.96                       & 3.25                        \\
                    & 4                     & 90.40                    & 213.38                      & 37.20                        & 98.20                    & 216.10                      & 33.97                        & 82.64                    & 193.07                      & 32.73                       \\ \hline
\multirow{3}{*}{FL} & 2                     & 12.40                    & 0.52                        & 0.29                         & 10.60                    & 0.51                        & 0.26                         & 9.43                     & 0.51                        & 0.30                        \\
                    & 3                     & 42.73                    & 15.32                       & 3.65                         & 42.91                    & 12.36                       & 3.51                         & 36.56                    & 12.17                       & 3.51                        \\
                    & 4                     & 87.16                    & 192.87                      & 38.87                        & 92.28                    & 163.37                      & 78.79                        & 77.73                    & 548.09                      & 48.33                       \\ \hline
\multirow{3}{*}{BL} & 2                     & 12.60                    & 0.53                        & 0.28                         & 10.74                    & 0.51                        & 0.30                         & 9.56                     & 0.61                        & 0.29                        \\
                    & 3                     & 43.04                    & 14.63                       & 4.10                         & 43.04                    & 15.15                       & 4.03                         & 36.86                    & 13.56                       & 3.76                        \\
                    & 4                     & 87.32                    & 225.83                      & 74.67                        & 92.33                    & 217.33                      & 71.23                        & 78.13                    & 190.33                      & 35.40                       \\ \hline
                    & \multicolumn{1}{l|}{} & \multicolumn{9}{c}{$c_o=200$, $c_r=50$, $c_e^1=150$, $c_e^2=100$}                                                                                                                                                                                                   \\ \hline
                    & \multicolumn{1}{l|}{} & \multicolumn{3}{c|}{Equal-segmented}                                                & \multicolumn{3}{c|}{High-priority}                                                    & \multicolumn{3}{c}{Low-priority}                                                     \\ \hline
                    &                       & \multicolumn{1}{c}{Cost} & \multicolumn{1}{c}{$DLPt(sec)$} & \multicolumn{1}{c|}{$RLPt(sec)$} & \multicolumn{1}{c}{Cost} & \multicolumn{1}{c}{$DLPt(sec)$} & \multicolumn{1}{c|}{$RLPt(sec)$} & \multicolumn{1}{c}{Cost} & \multicolumn{1}{c}{$DLPt(sec)$} & \multicolumn{1}{c}{$RLPt(sec)$} \\ \hline
\multirow{3}{*}{EL} & 2                     & 8.56                     & 0.74                        & 0.31                         & 10.54                    & 0.62                        & 0.35                         & 5.32                     & 0.60                        & 0.44                        \\
                    & 3                     & 24.96                    & 16.36                       & 3.76                         & 37.11                    & 14.54                       & 3.85                         & 17.61                    & 14.34                       & 3.65                        \\
                    & 4                     & 48.37                    & 306.64                      & 54.29                        & 76.64                    & 376.21                      & 43.29                        & 33.52                    & 185.91                      & 33.99                       \\ \hline
\multirow{3}{*}{FL} & 2                     & 5.85                     & 0.55                        & 0.31                         & 6.97                     & 0.59                        & 0.36                         & 3.74                     & 0.57                        & 0.46                        \\
                    & 3                     & 21.86                    & 13.41                       & 4.65                         & 32.3                     & 13.65                       & 4.28                         & 13.74                    & 12.89                       & 4.21                        \\
                    & 4                     & 45.71                    & 187.52                      & 41.15                        & 71.00                    & 174.55                      & 38.36                        & 28.79                    & 588.74                      & 34.15                       \\ \hline
\multirow{3}{*}{BL} & 2                     & 6.79                     & 0.59                        & 0.34                         & 7.66                     & 0.61                        & 7.66                         & 4.12                     & 0.64                        & 0.34                        \\
                    & 3                     & 23.62                    & 13.80                       & 4.71                         & 32.15                    & 13.61                       & 32.15                        & 15.17                    & 14.03                       & 4.42                        \\
                    & 4                     & 48.75                    & 220.22                      & 72.95                        & 70.22                    & 305.77                      & 70.22                        & 31.71                    & 170.47                      & 34.23                       \\ \hline \hline
\end{tabular}}
\label{Table::5}
\end{table}

\begin{table}[H]
\setcellgapes{3pt}\makegapedcells
\caption{Computational requirements for $K=2$, $M=2$.}
\centering
\resizebox{0.8\textwidth}{!}{
\begin{tabular}{r|rrrr|rr|rr}
\hline \hline
\multicolumn{1}{l|}{} & \multicolumn{4}{c|}{$c_o=200$, $c_r=150$, $c_e^1=300$, $c_e^2=250$}                                                                                                                                  & \multicolumn{4}{c}{$c_o=200$, $c_r=50$, $c_e^1=150$, $c_e^2=100$}                                                                                                              \\ \hline
\multicolumn{1}{l|}{} & \multicolumn{2}{c}{DLP}                                                                                     & \multicolumn{2}{c|}{RLP}                                                              & \multicolumn{2}{c|}{DLP}                                                               & \multicolumn{2}{c}{RLP}                                                                \\ \hline
$A$                   & \makecell{Number of \\actions} & \multicolumn{1}{r|}{$Mem$}  & \makecell{Number of \\actions} & $Mem$ & \makecell{Number of \\actions} & $Mem$  & \makecell{Number of \\actions} & $Mem$  \\ \hline
2                     & 1896                                                                          & \multicolumn{1}{r|}{4.08}   & 728                                                                           & 1.69  & 1896                                                                          & 4.08   & 864                                                                           & 1.97   \\
3                     & 12210                                                                         & \multicolumn{1}{r|}{74.93}  & 2952                                                                          & 18.63 & 12210                                                                         & 74.93  & 3262                                                                          & 20.52  \\
4                     & 53290                                                                         & \multicolumn{1}{r|}{778.97} & 8640                                                                          & 128.1 & 53290                                                                         & 778.97 & 9238                                                                          & 136.82 \\ \hline \hline
\end{tabular}}
\label{Table::6}
\end{table}

\end{document}